


\documentstyle{amsppt}

\define\r{{\Bbb R}}

\define\NS{\text{NS}}
\define\Exc{\text{Exc}}
\define\codim{\text{codim}~}
\define\rk{\text{rk}~}
\define\M{{\Cal M}}

\define\E{{\Cal E}}
\define\La{{\Cal L}}

\define\Ha{{\Cal H}}
\define\T{{\Cal T}}
\define\N{{\Cal N}}
\define\Pa{{\Cal P}}

\define\nem{\overline {NE}}

\magnification1200

\topmatter

\title
Algebraic surfaces with log-terminal singularities and $nef$
anticanonical class and reflection groups in Lobachevsky spaces.\\
I. (Basics of the diagram method)
\endtitle

\author
Viacheslav V. Nikulin \footnote{Partially supported by
Grant of Russian Fund of Fundamental Research;
Grant of AMS; and Grant of ISF M16000.}
\endauthor

\address
Steklov Mathematical Institute,
ul. Vavilova 42, Moscow 117966, GSP-1, Russia.
\endaddress

\email
slava\@nikulin.mian.su
\endemail

\date
May 22, 1994
\enddate

\abstract
In our preprint: "Algebraic surfaces with log terminal singularities
and $nef$ anticanonical class and reflection groups in Lobachevsky
spaces", Preprint Max-Planck-Institut f\"ur Mathematik, Bonn,
(1989) MPI/89-28 (Russian), we had extended the diagram method to algebraic
surfaces with $nef$ anticanonical class (before, it
was applied to Del Pezzo surfaces). The main problem here is
that Mori polyhedron may not be finite polyhedral in this case.
In the \S 4 of this preprint we had shown that this method
does work for algebraic surfaces with $nef$ anticanonical
class and log terminal singularities.

Recently, using in particular these our results, V.A. Alexeev
(Boundedness and $K^2$ for log surfaces, Preprint (1994),
alg-geom/942007), obtained strong and important results
for surfaces of general type with semi-log canonical
singularities. The reason is that Del Pezzo surfaces
(the diagram method is especially effective for them)
and surfaces of general type "are connected" through surfaces
with numerically zero canonical class. This was the main reason that
we decided to make this English translation of the first part
(\S 1---\S 3)
of our preprint above with some extensions (where we give
omitted proofs) and minor corrections.

Thus, this part contains basics of the diagram
method for generalized reflection groups in Lobachevsky spaces and
for surfaces with $nef$ anticanonical class.
\endabstract

\endtopmatter

\rightheadtext{Reflection Groups and Algebraic Surfaces}

\document

\vskip10pt

\noindent
Introduction

\noindent
\S 1 Convex polyhedra with acute angles in Lobachevsky spaces and
the diagram method

\noindent
1.1. Klein model of Lobachevsky space

\noindent
1.2. Polyhedra in Lobachevsky space

\noindent
1.3. Polyhedra with acute angles in Lobachevsky space

\noindent
1.4. Diagram method

\noindent
1.5. The proof of Theorem 1.4.0

\noindent
\S 2 Discrete reflection groups in Lobachevsky spaces and the diagram method

\noindent
\S 3 Diagram method and  algebraic surfaces

\noindent
\S 4 Projective algebraic surfaces with log terminal singularities and
$nef$ anticanonical class

\noindent
4.1 Two-dimensional log terminal singularities

\noindent
4.2 Projective algebraic surfaces with log terminal singularities and
$nef$ anticanonical class

\noindent
4.3 Arithmetic nature of the set of exceptional curves on Del Pezzo
surfaces with log terminal singularities and $nef$ anticanonical class

\noindent
4.4 Log terminal almost Del Pezzo and Del Pezzo systems
of exceptional curves and vectors

\noindent
4.5 Types and invariants of surfaces with $nef$ anticanonical class and
log terminal singularities

\noindent
4.6 Almost Del Pezzo surfaces with log terminal singularities

\noindent
4.7 Surfaces with $nef$ anticanonical class and log terminal singularities

\vskip20pt

\head
0. Introduction
\endhead

In fact, this is the English translation
with some extensions (where we give proofs which were
omitted) and minor corrections of the
first part \S 1---\S 3 of our preprint \cite{33}. This part
contains basics of the diagram method for generalized reflection
groups in Lobachevsky spaces and for algebraic surfaces with $nef$
anticanonical class.
We hope to translate the rest part of the preprint \cite{33} later.

In the series of papers of the author \cite{10}, \'E.B. Vinberg \cite{6},
A.G. Khovanskii \cite{15} and M.N. Prokhorov \cite{14}, it was proved that
discrete reflection groups in Lobachevsky spaces with the
fundamental polyhedron of
finite volume exist in dimension $<\ 996$ only.
The proof was based on the investigation of geometrical properties of
the fundamental polyhedra $\M$ of the groups. They are closed convex
polyhedra characterized by the property that any facial angle (the
angle between faces of highest dimension) has the
form $\pi /n,\ n\in \Bbb N,\ n>1.$ In \cite{11}, \cite{12} we remarked that
in the above series of papers a general result about polyhedra with acute
angles in Lobachevsky space was proved. We formulate this result now.

Let $\M$ be a convex locally-finite
polyhedron in a Lobachevsky space $\La $
(defined by a hyperbolic form $\Phi$ over $\r$
of signature $(1,n),\ n=\dim \La$)
and $P(\M)\subset \Phi$ be the set of vectors orthogonal to
faces of $\M$ of highest dimension and directed outside (one
vector for each face). Let $\Gamma (P(\M ))=(v_i\cdot v_j),
\ v_i,v_j \in P(\M)$ be the {\it Gram matrix} of
$\M$. We can also correspond to $P(\M)$ (and subsets of $P(\M)$ as well)
the {\it Gram graph} $\Gamma (P(\M ))$
with the set of vertices $P(\M)$ and edgers $v_iv_j$ if
$v_i\cdot v_j\not=0$ with
the weight $v_i\cdot v_j$. Evidently, the
Gram graph is equivalent to the Gram matrix, and we denote them by the
same letter. Using the Gram graph, we can consider paths,
length of a path, and so on.
A subset $\E\subset P(\M )$ is called {\it elliptic} if its Gram matrix
is negative definite. A subset $L \subset P(\M )$ is called {\it Lanner}
if its Gram matrix is hyperbolic, but the Gram matrix of any
proper subset of $L$ is not hyperbolic (or $L$ is a minimal hyperbolic
subset).

A polyhedron $\M$ {\it has acute angles} if $v_i\cdot v_j\ge 0$
for any $v_i\not=v_j \in P(\M)$. A polyhedron $\M$ is called
{\it closed} if it is a convex envelope of a finite set of points
(with some of them at infinity) in
$\La$. The polyhedron $\M$ is called {\it non-degenerate} if it is not
contained in a hyperplane.
In fact, in papers \cite{10}, \cite{6}, \cite{15}, \cite{14}, the
following general result about polyhedra with acute angles in a
Lobachevsky
space was obtained.

\proclaim{Lemma 1} Let $\M$ be a
closed non-degenerate polyhedron with acute angles in
a Lobachev\-sky space $\La$ and suppose there exist some constants $d, C_1,
C_2$ such that

(a) $\text{diam\ } \Gamma (L) \le d$ for any Lanner subset
$L \subset P(\M )$;

(b)
$$
\#\{ \{ v_1,v_2\} \subset \E \ \mid \ 1 \le \rho (v_1,v_2)\le d \} \le
C_1\# \E
$$
and
$$
\#\{ \{ v_1,v_2\} \subset \E \ \mid \ d+1\le \rho (v_1,v_2)\le 2d+1\} \le
C_2\# \E
$$
for any elliptic subset $\E \subset P(\M )$ and the distance in
$\Gamma (\E )$.
Then
$$
\dim \La < 96(C_1+C_2/3)+68.
$$
\endproclaim

Here we have formulated slightly stronger result than was formulated in
\cite{11}, \cite{12} (we consider $\text{diam\ }\Gamma( L)$ instead of
$\sharp L$.

One of the key points of this paper is
to generalize Lemma 1
for some convex locally-finite polyhedra with acute angles in
Lobachevsky spaces which
are not necessarily closed
(in particular, they may have infinite volume and may have
infinite number of faces). To do this, we need several definitions
(see Sect. 1.2).

A convex polyhedron $\M$ is called {\it elliptic}
if it is closed (i.e. it is the convex hull of a finite set of points)
and is non-degenerate. The same definition works for Euclidean space too.

A convex polyhedron $\M$ in a Lobachevsky space $\La$ is called
{\it parabolic} relative to a point $P$ at infinity if for any elliptic
polyhedron $K$ on the horosphere with the center $P$ the polyhedron
$C_K\cap \M$ is elliptic if it is non-degenerate.
Here $C_K$ is the cone with the vertex $P$ over $K$.

A convex polyhedron $\M$ in
a Lobachevsky space $\La$ is called {\it hyperbolic}
relative to a subspace
$\T\subset \La$ with
$0 < \dim \T<\dim \La$ if it is finite in infinite points
and for any elliptic compact polyhedron
$K \subset \T$ the polyhedron $C_K \cap \M$ is elliptic. Here $C_K$ is
the union of the lines $AB, A\in K$ and $AB\perp \T$ (the cylinder with
the base $K$).

We prove (Lemmas 1.2.1 and 1.2.2) that if $\M$ is convex parabolic,
then there exists an elliptic face $\gamma \subset \M$ of codimension $1$
with $P\notin [\gamma]$. If $\M$ is convex hyperbolic, then there exists an
elliptic face $\gamma \subset M$ of $\codim\gamma \le \dim \T$,
and $\gamma$ is cut out by faces hyperplanes of $\M$ which do not
contain $\T$. If $\M$ has acute angles then $\gamma$ also has acute
angles, and we can apply Lemma 1 to $\gamma$. By this way, we can apply
Lemma 1 (and its modifications)
to elliptic, parabolic and hyperbolic convex polyhedra
with acute angles in Lobachevsky spaces. See Sect. 1.4, Theorems
1.4.0 --- 1.4.3, 1.4.2', 1.4.3' for
exact statements (also see Theorem 1.5.1.1. for the compact case).
These statements constitute the diagram method for polyhedra with acute
angles in Lobachevsky spaces which are not necessarily closed.
Below we show that these statements (diagram method)
are suitable for applications to
appropriate classes of reflection groups in Lobachevsky spaces and
algebraic surfaces.

In \S 2 we apply the diagram method
to so called generalized (hyperbolic)
crystallographic reflection groups $G$ in Lobachevsky spaces
(probably, at first they were introduced in \cite{26}).
Let $G$ be a discrete group in a Lobachevsky space $\La$ with a
fundamental domain of finite volume. We say that $G$ is a
{\it generalized crystallographic reflection group} if for its
reflection subgroup $W$ (generated by all reflections in
hyperplanes which belong to $G$) the quotient group $G/W$, considered as a
group of symmetries of some fundamental polyhedron $\M$ of $W$,
contains a subgroup $A^\prime$ of a finite index which leaves invariant
a proper subspace $\La^\prime$ of $\La$. We show (see Proposition 2.1)
that in this definition we can suppose that $\La^\prime$ is
generated by $\La^\prime \cap \M$. We prove (Proposition 2.2) that
in this case $\M$ is elliptic if $\La^\prime$
is a finite point, $\M$ is parabolic relative to $\La^\prime$ if
$\La^\prime$ is an infinite point, and $\M$ is
hyperbolic relative to
$\La^\prime$ if $\La^\prime$ is a finite subspace of $\La$.
Thus, we can apply the diagram method to
the fundamental polyhedron $\M$ (since it has acute angles). We get

\proclaim{Theorem 2.3} Let $G$ be a generalized crystallographic
reflection group in a Lobachevsky space $\La$. Then

(a) $\dim \La < 996$ if $\La^\prime$ is a finite point
(this is the result of A.G. Khovanskii \cite{15} and
M.N. Prokhorov \cite{14});

(b) $\dim \La < 997$ if $\La^\prime$ is an infinite point;

(c) $\dim \La < 996 + \dim \La^\prime $ if $\La ^\prime$ is a finite
subspace of $\La$.
\endproclaim

We mention that if $\M$ is bounded,
\'E.B. Vinberg \cite{6} had proven that $\dim \La < 30$.
Before, the author had proven \cite{10} that $\dim \La < 10$
if $G$ is an arithmetic reflection group and the degree
of the ground field is sufficiently large. Thus, Theorem 2.3 generalizes
these results of the author, \'E.B. Vinberg, A.G. Khovanskii and
M.N. Prokhorov to wider class of groups. One can find non-trivial
examples of generalized crystallographic reflection groups in
\cite{9}.

We don't know whether there exists an absolute estimate of $\dim \La$
(which does not depend from $\dim \La^\prime$) for
the case (c).

\smallpagebreak

In \S 3 we show that the diagram method above may be applied to
appropriate classes of algebraic surfaces.

In the first place, we consider non-singular projective
algebraic surfaces $Y$ over an algebraically closed field.
Let $\NS(Y)$ be the Neron-Severi lattice of $Y$ (generated by algebraic
curves by the numerical equivalence). By Hodge index
theorem, the lattice $\NS (Y)$ is hyperbolic
(i.e. it is non-degenerate and has exactly one
positive square) and defines the Lobachevsky
space
$\La (Y)= V^+ (Y)/\r^+$  where $V(Y)=\{x\in \NS (Y)\otimes \r \mid
x^2>0\}$ and $V^+(Y)$ is the half-cone containing
the class of a hyperplane section of $Y$.

We recall that a curve $C$ of $Y$ is called {\it exceptional} if
$C$ is irreducible and $C^2<0$. Let $\Exc(Y)$ be the
set of all exceptional curves of $Y$. Any
exceptional curve $F$ defines a half-space
$\Ha_F^+=\{ \r^+x \in \La (Y) \mid x\cdot F\ge 0\}$
bounded by the hyperplane
$\Ha_F=\{ \r^+x \in \La (Y) \mid x\cdot F = 0\}$
in $\La (Y)$. The set $\Exc (Y)$ defines a convex polyhedron
$$
\M(Y)=\bigcap_{F \in \Exc(Y)}{\Ha_F^+}
$$
in $\La (Y)$ with $P(\M (Y))=\Exc (Y)$.
The $\M(Y)$ has acute angles
since $F \cdot F^\prime \ge 0$ if $F\not= F^\prime
\in \Exc (Y)$.

 From the classical results of Mori \cite{23} for the case of
surfaces, we deduce

\proclaim{Proposition 3.1} Let $Y$ be a non-singular projective
algebraic surface
over an algebraically closed field. Assume that the anticanonical
class $-K_Y$ is pseudoeffective (i.e. $-K_Y$ belongs to the Mori cone
$\nem (Y)$) and $-K_Y \not\equiv 0$. Let
$$
-K_Y \equiv P+\sum_{i=1}^n {\alpha_i F_i}
$$
be Zariski decomposition of $-K_Y$
(here $P$ is $nef$; $P\cdot F_i=0$ and $\alpha_i>0$ for all $i$;
curves $F_i,\ i=1,2,...,n$, have a negative definite Gram matrix).
Then $\M (Y)$ is elliptic if $P^2>0$, $\M(Y)$
is parabolic relative to $\Bbb R^+P$ if $P^2=0$ and
$P\not\equiv 0$, and $\M(Y)$ is hyperbolic relative to the
subspace $\T =\cap \Ha_{F_i} \subset \La $ (intersection of hyperplanes
$\Ha_{F_i}$ orthogonal to $F_i,\ i=1,...,n$) of $\codim \T=n$
if $P\equiv 0$.
\endproclaim

Thus, we can apply the diagram method above
to non-singular projective algebraic surfaces
with pseudo-effective anti-canonical class $-K_Y$. See exact
formulations in Theorems 3.4---3.6 and 3.5', 3.6'.

These applications
are especially important for projective algebraic surfaces with
normal (or isolated ) singularities and
$nef$ anti-canonical class.

Let $Z$ be a normal projective algebraic surface,
and $\sigma : Y\to Z$ the minimal resolution of
singularities of $Z$. Let the anticanonical class
$-K_Z$ is $nef$. Then Zariski decomposition of $-K_Y$ has the
form
$$
-K_Y\equiv -\sigma ^\ast K_Z+\sum {\alpha_jF_j},
\tag0.1
$$
where $F_j$ are the components of the exceptional divisor of $\sigma$.
Here all $\alpha_j\ge 0$ and $\alpha_j=0$ iff $\sigma (F_j)$ is
a double rational point (Du Val singularity). Thus, to get Zariski
decomposition, we should just take away these zero summands in (0.1).

By Proposition 3.1, $\M(Y)$ is elliptic if $(K_Z)^2>0, \M(Y)$ is parabolic
if $(K_Z)^2=0$ and $K_Z\not\equiv 0$, $\M (Y)$ is hyperbolic if
$K_Z\equiv 0$. Hence we can apply
the diagram method to $\M(Y)$. Thus, using the diagram method
we get an estimate for $\dim \NS(Y)$ if $K_Z\not\equiv 0$.
And we get an estimate for
the number $n$ of curves $F_j$ with positive $\alpha_j$ in the
formula (0.1) if
$K_Z \equiv 0$. At any case, we get the non-trivial
restriction on singularities of $Z$.

It is an interesting problem to find classes of singularities of $Z$
for which the results above (the diagram method)
give non-trivial and interesting applications.
This paper contains only
basics of the diagram method, and we don't discuss these
applications. Here in Introduction
we can only give references where
one can find these applications.

In our papers \cite{11}, \cite{12}, \cite{13} and of V.A. Alexeev
\cite{2} the diagram method was applied very effectively to
Del Pezzo surfaces (when $-K_Z$ is numerically ample)
with log terminal singularities (for surfaces over $\Bbb C$
"log terminal" is the same as "quotient").
For this case, the polyhedron $\M(Y)$ is
elliptic, and one can apply Lemma 1 directly.

For projective algebraic surfaces
with log-terminal singularities and $nef$ anticanonical class,
the diagram method presented here was first applied in the rest part
(\S 4) of the preprint \cite{33}. We hope to translate and
prepare for publishing this part later.

Recently, V.A. Alexeev \cite{37} (see also \cite{36})
found very interesting and important applications of the diagram method
for algebraic surfaces with $nef$ anticanonical
class which we present here,
to algebraic surfaces of general type.
Probably, the reason is that surfaces
of positive canonical class "are connected" with surfaces
of negative canonical class through surfaces with numerically
zero canonical class. Actually, this results of
V.A. Alexeev were the reason that
we decided to make this English
translation of our preprint \cite{33}.

Some generalization of the diagram method to $3$-folds one can
find in \cite{34} and \cite{35}.

I am grateful to Professor V.A. Alexeev
for his interest to this translation.

\newpage

\head
\S 1. Convex polyhedra with acute angles
in Lobachevsky spaces and the diagram method
\endhead

Hear we obtain some general results about polyhedra in Lobachevsky space.
They give the basics of diagram method we will later apply to reflection
groups in Lobachevsky spaces and to algebraic surfaces.

\subhead
1.1. Klein model of Lobachevsky space
\endsubhead

Let $\Phi$ be a hyperbolic linear space (i.e. $\Phi$  is a
$\Bbb R$-linear space equipped with a non-degenerate
symmetric $\Bbb R$-bilinear form
with inertia indexes $(1, n)$, where $n+1=\dim \Phi$).
For $x,y \in \Phi$,
we denote as $x\cdot y$
the value at the pair $(x,y)$ of the corresponding to $\Phi$
symmetric bilinear form. One can associate to $\Phi$ an open cone
$$
V=\{x \in \Phi \ \mid \ x^2>0\}.
$$
Since the bilinear form of $\Phi$ is hyperbolic,
the $V$ is the disjoint union of two open convex half-cones.
We denote by $V^+$ one of them. Thus, $V=V^+ \cup -V^+$.
The corresponding to $\Phi$ Lobachevsky space is the set of
rays
$$
\La =V^+/\r^+ = \{ \r^+ x\ \mid \ x\in V^+ \},
$$
where $\r^+$ denote the set of positive real numbers.
The distance $\rho$ in $\La$ is defined by the formula
$$
\cosh \rho (\r^+x, \r^+y)=(x\cdot y)/(x^2y^2)^{1/2}.
$$
Then the curvature of $\La$ is equal to (-1).

Let $O(\Phi)$ be the group of automorphisms of the hyperbolic linear
space $\Phi$. Then the group
of motions of $\La$ is identified with the subgroup
$O_+(\Phi)\subset O(\Phi)$ of index $2$
of automorphisms which preserve the half-cone $V^+$.

A half-space of $\La$ is the set
$$
\Ha_\delta^+ = \{\r^+x \in \La\ \mid \ x\cdot \delta \ge 0\} ,
$$
where $\delta \in \Phi$ and $\delta^2<0$. The half-space $\Ha_\delta^+$
is bounded by the hyperplane
$$
\Ha_\delta = \{\r^+x \in \La\ \mid \ x\cdot \delta = 0 \}.
$$
The element $\delta\in \Phi$ is defined by the half-space $\Ha_\delta^+$
(respectively, the hyperplane $\Ha_\delta$)
up to multiplying by elements of $\r^+$ (respectively, by elements of
the set $\r^\ast$ of non-zero real numbers). The $\delta$ is called
orthogonal to the half-space $\Ha_\delta^+$ (respectively,
the hyperplane $\Ha_\delta$).

Assume that intersection
$\Ha_{\delta_1}^+\cap \Ha_{\delta_2}^+$ of two
half spaces $\Ha_{\delta_1}^+,\  \Ha_{\delta_2}^+$ orthogonal to
$\delta_1,\ \delta_2 \in \Phi$ contains a non-empty open subset of $\La$.
Then we have two cases:

a) $\Ha_{\delta_1}^+\cap \Ha_{\delta_2}^+$ is the angle of the
value $\phi$ where
$$
\cos \phi =(\delta_1\cdot \delta_2)/(\delta_1^2 \delta_2^2)^{1/2}
\tag1.1
$$
if $-1\le (\delta_1\cdot \delta_2)/(\delta_1^2 \delta_2^2)^{1/2}\le 1.$

b) Intersection of the hyperplanes $\Ha_{\delta_1}, \Ha_{\delta_2}$
is empty (they are hyperparallel),
and the distance $\rho$ between them is defined by
$$
\cosh \rho =(\delta_1\cdot \delta_2)/(\delta_1^2 \delta_2^2)^{1/2}
\tag1.2
$$
if $1\le (\delta_1\cdot \delta_2)/(\delta_1^2 \delta_2^2)^{1/2}.$

As usual, we complete $\La$ by the set of infinite points
$\r^+c$ where $c\in \Phi,\ c^2=0, c\cdot V^+>0$. Then $\La$ with its
infinite points is the closed ball
$$
\overline\La=(\overline V^+-\{ 0\} )/\r^+,
$$
where $\overline V^+$ is the closure of $V^+$
in the linear space $\Phi$). The boundary
$\La_\infty = \overline\La - \La$
is the sphere of the dimension $n-1$. It is
called the $\infty$-sphere. We also add infinite points
to half-spaces and hyperplanes considering their closure in
$\overline \La$. Thus,
$$
\overline \Ha_\delta^+=
\{\r^+x \in \overline \La\ \mid \ x\cdot \delta\ge 0\},
$$
and
$$
\overline \Ha_\delta =
\{\r^+x \in \overline \La\ \mid \ x\cdot \delta = 0\}
$$
are the half-space $\Ha_\delta^+$ and the hyperplane
$\Ha_\delta$ with their infinite points respectively.

\subhead
1.2. Polyhedra in Lobachevsky space
\endsubhead

A non-degenerate convex locally finite polyhedron $\M$
in a Lobachevsky space $\La$ is an intersection of a set of half-spaces:
$$
\M=\bigcap_{\delta \in P(\M )}{\overline \Ha_\delta^+}.
$$
The $\M$ is called {\it non-degenerate} if $\M$ contains a non-empty open
subset of $\La$. The $\M$ is called {\it locally finite} if for any
$X\in \La$, there exists its open neighborhood $U \subset \La$
such that $\M\cap U$ is intersection with $U$
of a finite set of half-spaces in $\La$.
Here $P(\M)\subset \Phi$ is a finite or countable subset of elements
with negative square. We assume that neither
two of these elements are proportional and
every half-space $\Ha_\delta^+,\
\delta \in P(\M)$ defines a face of $\M$ of the highest dimension
$\dim \La-1$ (it is clear what this means, since $\M$ is
locally finite). Then, the polyhedron $\M$ defines the subset
$P(\M )\subset \Phi$ uniquely up to multiplying by $\r^+$ elements
of $P(\M)$.

Below, we are only considering locally-finite and convex polyhedra.

We recall that a face $\gamma$ of a polyhedron $\M$ is intersection
$$
\gamma = \M \cap (\bigcap_{\delta \in S}{\Ha_\delta})
$$
where $S$ is a non-empty subset of $P(\M)$. One can consider $\gamma$
as a polyhedron is the subspace
$[\gamma ]$ of $\La$ generated by $\gamma$ if $\gamma$ is finite
(i.e. contains a finite point of $\La$).
Otherwise, $\gamma$ is either infinite (is a point of $\La_\infty$)
or is empty.

\subsubhead
1.2.1. Elliptic polyhedra
\endsubsubhead

A convex non-degenerate polyhedron $\M$ in Lobachevsky space
$\La$ is called {\it elliptic} (equivalently {\it closed})
if $P(\M)$ is finite (in particular, $\M$ is finite)
and $\M\cap \La_\infty$ is finite.
Equivalently, $\M$ is the convex envelope of a finite set of points of
$\overline \La$ which is not contained in a hyperplane. The last
definition is valid for Euclidean space too.

A face $\gamma$ of $\M$ is called {\it elliptic} if $\gamma$ is an
elliptic polyhedron in $[\gamma]$.
We say that a polyhedron $\M$ in $\La$ is
{\it elliptic} (equivalently, {\it closed}) {\it in
a neighborhood of its face $\gamma$} if there
exists a neighborhood $\overline {\gamma}\subset U$ in
$\overline\La$
and an elliptic polyhedron $\M^\prime$ in $\La$ such that
$\M \cap U=\M^\prime \cap U$. Obviously, the face $\gamma$
is an elliptic polyhedron in $[\gamma]$ and the set
$$
P(\gamma)=\{ \delta \in P(\M)
\mid \overline \Ha_\delta \cap \overline \gamma \not=\emptyset \}
$$
is finite if $\M$ is elliptic
in a neighborhood of its face $\gamma$.

\subsubhead
1.2.2. Parabolic polyhedra
\endsubsubhead

Let $P=\r^+c\in \La_\infty$ ($c\in \Phi$ and $c^2=0$ but $c\not=0$)
is a point at infinity of $\La$. We remind that the horosphere
$\E_P$ with the center $P$ is the set of all lines in $\La$
containing the $P$.
The line $l=P\r^+h\in \E_P,\ \r^+h\in \overline \La,$
is the set
$\overline l = \{ \r^+(tc+h)\ \mid \ t\in \r \ \text{and} \ (tc+h)^2>0$.
We fix a constant $R>0$. Then there exists a unique $\r^+h\in l$
such that $h\cdot c=R$ and $h^2=1$. Let $l_1,l_2\in \E_P$ and
$h_1, h_2$ the corresponding elements respectively
we had defined above. Let
$$
\rho(l_1,l_2)=\sqrt{-(h_1-h_2)^2}.
\tag2.1
$$
The horosphere $\E_P$ with this distance is an affine Euclidean space.
If one changes the constant $R$, the distance $\rho$ will be multiplied
by a constant. The set
$$
\E_{P,R}=\{ \r^+h \in \La\ \mid \ h\cdot c=R\ \text{and} \ h^2=1 \}
\cup \{ P \}
$$
is a sphere in $\overline \La$ touching $\La_\infty$ at the $P$.
Besides, the sphere $\E_{P,R}$
is orthogonal to a line $l\in \E_P$ at the
corresponding to the $l$ point $\r^+h, h\in \E_{P,R}$. The distance
of $\La$ induces an Euclidean distance in $\E_{P.R}$ which is similar
to the distance (2.1). The set $\E_{P,R}$ is identified with
$\E_P$ and called horosphere too.

Let $K\subset \E_P$. The set
$$
C_K=\bigcup_{l\in K}{\overline l}
$$
is called the cone with the vertex $P$ and the base $K$. Evidently,
the cone $C_K$ is a subspace of the Lobachevsky space $\La$ if the
$K$ is a subspace of the Euclidean space $\E_P$. The dimension
$\dim C_K = \dim K +1$.

A non-degenerate locally finite polyhedron $\M$ in $\La$ is called
{\it parabolic} (relative to the point $P\in \La_\infty$), if
the conditions 1) and 2) below are valid:

1) $\M$ is finite at the point $P$, i.e. the set
$\{\delta \in P(\M) \mid c\cdot \delta =0\}$ is finite.

2) For any elliptic polyhedron ${\Cal N}\subset \E_P$ (i.e. $\Cal N$ is
a convex envelope of a finite set of points in $\E_P$), the
polyhedron $\M\cap C_{\Cal N}$ is elliptic if it is non-degenerate.

By the definition, a parabolic polyhedron $\M$ is elliptic
if and only if either $P\notin \M$ or
$$
\M_P = \bigcap\limits_{\delta \in P(\M),\ \delta\cdot c=0}
{\E_P \cap \Ha_\delta^+}
\tag2.2
$$
is an elliptic polyhedron in $\E_P$.

We need

\proclaim{Lemma 1.2.1} Let $\M$ be a non-degenerate
locally-finite convex polyhedron in a
Loba\-chevsky space $\La$ and
$\M$ is parabolic relative to a point $P\in \La_\infty$.
Then there exists a face $\gamma$ of $\M$ such that
$\dim \gamma =\dim \La -1$,
$\overline {[\gamma]}$ does not contain $P$ and $\M$ is
elliptic in a neighborhood of $\gamma$ (in particular,
$\gamma$ is elliptic in the hyperplane $[\gamma]$).
\endproclaim

\demo{Proof}
Assume that $P\notin \M$.
Then, there exists $\delta \in P(\M)$ such that
$P \notin \overline \Ha_\delta^+$. The set
$D=\{ PX\in \E_P \ \mid \ X\in \overline \Ha_\delta^+\}$
is a ball in $\E_P$. There exists an elliptic polyhedron
${\Cal N}\subset \E_P$ which contains the $D$. Then
$\M=\M \cap C_{\Cal N}$ is elliptic, and the face $\gamma$ of $\M$
we are looking for does evidently exist.

Assume that $P\in \M$. Let $Q$ be a point inside of $\M$. There
exists an elliptic polyhedron
${\Cal N}\subset \E_P$ such that ${\Cal N}\subset \M_P$ (see (2.2))
and $PQ\in {\Cal N}$ is an internal point of ${\Cal N}$. A polyhedron
$C_{\Cal N}\cap \M$ is elliptic, since $\M$ is parabolic.
By the construction, the $C_{\Cal N}\cap \M$ and the cone $C_N$ coincide
at the point $P$. It follows that there exists a codimension one
face $\beta$ of $C_{\Cal N}\cap \M$ such that
$P\notin \overline{[\beta]}$. By construction, the $\overline {[\beta]}$
is the hyperplane of a codimension one face $\gamma$ of $\M$,
$\beta \subset \gamma$ and $[\gamma]=[\beta]=\Ha_\delta$ for
$\delta \in P(\M),\ \delta \cdot c\not=0$ (here $P=\r^+c$).
The set
$$
C_\gamma =\{l\in \E_P \ \mid \ l\cap \overline \Ha_\delta \not=\emptyset\}
$$
is a closed ball of the horosphere $\E_P$. There exists an elliptic
polyhedron $T\subset \E_P$ such that the ball $C_\gamma$ is contained
in interior of the $T$. The polyhedra $\M$ and $\M\cap C_T$ coincide in
a neighborhood of $\overline \Ha_\delta$. Since $\M\cap C_T$ is
elliptic (because the $\M$ is parabolic),
the polyhedron $\M$ is elliptic in the neighborhood of
$\overline \Ha_\delta$.
$\blacksquare$
\enddemo

\subsubhead
1.2.3. Hyperbolic polyhedra
\endsubsubhead

Let $\T$ be a Lobachevsky subspace (i.e. intersection of a finite set of
hyperplanes) of a Lobachevsky space $\La$, and
$1\le \dim \T \le \dim \La -1$. Let $K \subset \T$. We denote
$$
C_K=\bigcup_{l\perp \T, \atop l\cap \T \in K}{l}
$$
where $l$ is a line of $\La$.
Thus, $C_K$ is the cylinder over the $K$.

A non-degenerate locally-finite polyhedron $\M$ in $\La$ is called
{\it hyperbolic} (relative to $\T$) if the conditions 1) and 2)
below hold:

1) $\M$ is finite in infinite points, i.e. the set
$\{\delta \in P(\M) \mid \delta\cdot c=0\}$ is finite for any point
$\r^+c \in \La_\infty$;

2) For any compact elliptic polyhedron $\N\subset \T$, the
polyhedron $C_\N\cap \M$ is elliptic if it is non-degenerate.

We need

\proclaim{Lemma 1.2.2} Let $\M$ be a non-degenerate locally-finite
convex polyhedron in a Lobache\-vsky space $\La$ and $\M$ is hyperbolic
relative to a subspace $\T \subset \La$ where
$0 < \dim \T< \dim \La$.
Then there exists  an elliptic face $\gamma$ of $\M$ such that
$\codim \gamma \le \dim \T$
and $\gamma$ is cut out by several hyperplanes
$\Ha_\delta,\ \delta \in P(\M),$ which don't contain the $\T$.
\endproclaim

\demo{Proof} We prove this by induction on $\dim \T$.

Let $\dim \T = 1$.
Assume that there exists a hyperplane $\Ha_\delta,\
\delta \in P(\M),$ such that either
$\overline \Ha_\delta \cap \overline \T =\emptyset$, or
$\Ha_\delta$ and $\T$ intersect in a finite point of $\La$. Let
$D\subset \T$ be the image of $\overline \Ha_\delta$
by the orthogonal projection on $\T$. Under the conditions above,
the $D$ is a compact interval.
Let $\N$ be an elliptic compact polyhedron
in $\T$ such that $D$ is contained inside of $\T$. By the construction,
the polyhedra $C_\N\cap \M$ and $\M$ coincide in a neighborhood of
$\overline \Ha_\delta$. Since $\M$ is hyperbolic, it follows that
$C_N\cap \M$ and $\M$ are both elliptic in the neighborhood of
$\overline \Ha_\delta$. Obviously, the codimension one face
$\gamma =\M\cap \Ha_\delta$ is the required one.
Now, assume that any hyperplane $\Ha_\delta,\ \delta\in P(\M),$ either
contains $\T$ or intersects $\T$ in infinite point.
In both cases, $\Ha_\delta$ contains an infinite point of $\T$.
Since $\M$ is finite at infinite points
(by the condition 1 above) and $\T$
contains exactly two infinite points, it follows that $P(\M)$ is finite.
Since $\M$ is hyperbolic relative to $\T$, one can easily see that
$\M$ is then elliptic. Obviously, then there exists
$\Ha_\delta,\ \delta \in P(\M),$ such that $\T\not\subset \Ha_\delta$.
The codimension one face $\gamma$ of $\M$ containing in the
$\Ha_\delta$ is required.

Let $\dim \T>1$. We consider several cases.

Assume that there exists a hyperplane $\Ha_\delta,\ \delta \in P(\M)$,
such that $\overline \T\cap \overline \Ha_\delta =\emptyset$.
Arguing like above, we find the required face
$\gamma\subset \Ha_\delta$ of $\M$ with
$\codim \gamma =1$.

Assume that there exists a hyperplane $\Ha_\delta,\ \delta\in P(\M)$,
such that $\overline \T \cap \overline \Ha_\delta$ is an infinite point
$P\in \La_\infty$. Let $\beta\subset \Ha_\delta$
be a face of $\M$ of codimension one. Let $A$ be an internal point of
$\beta$. There exists a half-space
$\Ha_e^+$ of $\La$ such that $P\notin \Ha_e^+$  but
$\ A \in \Ha_e^+$ and the
boundary $\Ha_e$ is orthogonal to $\T$. Let $D$ be the image by the
orthogonal projection onto $\T$ of the set
$\overline \Ha_e^+ \cap \overline \Ha_\delta$.
The $D\subset \T$ is compact, and there exists an elliptic compact
polyhedron $\N\subset \T$ such that
$D$ is contained inside of $\N$. The polyhedron
$C_\N\cap \M$ is elliptic, since $\M$ is hyperbolic. Then the polyhedron
$\Ha_e^+\cap C_\N\cap \M$ is elliptic either. By the construction,
the $\Ha_e^+\cap C_\N\cap \M$ and $\Ha_e^+\cap \M$ coincide in a
neighborhood of $\overline \Ha_e^+\cap \overline \Ha_\delta$. The
$\beta\cap \overline \Ha_e^+$ is  the face of $\Ha_e^+\cap \M$ of
codimension one (since the point $A$ is contained inside of
$\beta$ and $\Ha_e^+$),
and this face is contained in $\overline \Ha_\delta$.
Thus, the polyhedron $\Ha_e^+\cap \M$ is elliptic in a neighborhood of
its face $\beta\cap\overline \Ha_e^+$. Moving the point $A$
inside $\beta$ a little,
we can suppose that the line $PA\subset \Ha_\delta$
contains internal points (for example, the point $A$) of
the elliptic polyhedron $\beta\cap \Ha_e^+$ and does not
contain its infinite points. Then $PA\cap \Ha_e^+\cap \beta =[R,S]$
where $R\not=S$ and $R\in [P,S]$. There exists a face $\gamma^\prime$
of $\Ha_e^+\cap \beta$ which intersects the line $PA$ at the point $S$
and $\codim \gamma^\prime =2$. In particular,
$P\notin \overline {[\gamma ^\prime]}$. The subspace $[\gamma ^\prime]$ is
different from $\Ha_e\cap \Ha_\delta$ because the last subspace
intersects the line $PM$ at the point of $[P,R]$. It follows that
$\gamma ^\prime \subset \gamma$ where $\gamma$ is a face of
the face $\beta$ of $\M$, and $\codim \gamma =2$.
Thus, $P\notin \overline{[\gamma^\prime]}=\overline{[\gamma]}$.
It follows that
the image of $\overline{[\gamma]}$ by the
orthogonal projection onto $\T$ is compact in $\T$. Since
$\M$ is hyperbolic, it follows like above that $\M$ is elliptic
in a neighborhood of $[\gamma]$.
Let $\gamma \subset \Ha_\alpha ,\ \alpha \in P(\M)$. Let us
suppose that
$\T\subset \Ha_\alpha$. Then $P\in \Ha_\alpha$.
Since $\codim \gamma =2$ and $P\notin \overline{[\gamma]}$,
it follows
that $\Ha_\alpha = \Ha_\delta$. We get a contradiction since
$\Ha_\delta$ does not contain $\T$.
Thus, $\gamma$ is the required face of $\M$.

Assume that there exists $\Ha_\delta,\  \delta \in P(\M)$
such that $\Ha_\delta$ intersects $\T$. Then $\T^\prime=
\T\cap \Ha_\delta \subset \Ha_\delta$ is a subspace of
$\dim \T^\prime =\dim \T -1\ge 1$. One can easily see that
the codimension one face $\M^\prime$ of
$\M$ containing in $\Ha_\delta$  is a hyperbolic
polyhedron relative to $\T^\prime$. By induction
hypothesis, there exists an elliptic face
$\gamma$ of the $\M^\prime$
such that
$\dim \Ha_\delta - \dim \gamma =
\dim \La - 1 - \dim \gamma \le \dim \T^\prime=\dim \T -1$
and $\gamma$ is intersection of some
codimension one faces of the polyhedron
$\M^\prime$ in $\Ha_\delta$ which don't contain $\T^\prime$.
It follows that $\gamma$ is an elliptic face of $\M$,
$\codim \gamma \le \dim \T$, and $[\gamma ]$ is intersection of
codimension one faces $\Ha_\alpha, \alpha \in P(\M),$ which
don't contain $\T$.

If $\M$ does not have a hyperplane $\Ha_\delta,\ \delta \in P(\M),$
which satisfies one of conditions above, then all hyperplanes
$\Ha_\delta,\ \delta \in P(\M),$ contain the subspace $\T$.
Obviously, then $\M$ is not hyperbolic relative to $\T$ since
the condition 2) above is not valid for any compact polyhedron
$\N \subset \T$.
$\blacksquare$
\enddemo

\subhead
1.3. Polyhedra with acute angles in Lobachevsky space
\endsubhead

We recall basic facts about polyhedra with acute angles in a
Lobachevsky space
(for example, one can find a very elementary exposition in \cite{7}).

Below we assume that all polyhedra we are considering are
non-degenerate, locally-finite and convex.

We say that a polyhedron $\M$ in a Lobachevsky space $\La$
{\it has acute angles}
if $\delta \cdot \delta^\prime \ge 0$ for any
$\delta \not= \delta ^\prime \in P(\M)$. By (1.1) and (1.2), facial
angles (angles between codimension one faces) of $\M$ are really
acute ($\le \pi /2$). The inverse statement is also valid
(see E.M. Andreev \cite{5}), but we don't need this,
and use the definition above. It is not difficult to see that
any finite face of a polyhedron with acute angles
is a polyhedron with acute angles too.

The most important invariant of $\M$ is {\it Gram matrix}
$(\delta_i\cdot \delta_j)$, $\delta_i, \delta_j \in P(\M)$. This is
a symmetric matrix with negative diagonal
and non-negative non-diagonal
elements. This matrix (i.e the corresponding symmetric bilinear form)
is not more than hyperbolic (i.e. it has
not more than one positive square),  and has rank
$\le \dim \La +1$ because $P(\M)$ is a subset of the
hyperbolic space $\Phi$ of $\dim \Phi= \dim \La +1$.

We also consider the
corresponding to $\M$ {\it Gram graph} $\Gamma (P(\M))$ or $\Gamma (\M)$.
Vertices of $\Gamma (\M)$ correspond to elements of $P(\M)$.
A vertex $\delta \in P(\M)$ has the weight $(-\delta ^2)$.
Two vertices $\delta \not=\delta^\prime$ are
joined by an arrow of the weight $\delta \cdot \delta ^\prime$ if
$\delta\cdot \delta^\prime >0$. Evidently, the Gram graph is equivalent
to the Gram matrix, and we identify them. Similarly, we consider
Gram matrix and graph for any subset of $P(\M)$.
Using Gram graph, we can consider connected components,
paths, length of paths, and so on.

A subset $\E\subset P(\M)$ is called {\it elliptic} if its Gram
matrix (or graph) is negative definite. These subsets are important
because they are in one to one correspondence with finite (i. e. having a
finite point of $\La$) faces of $\M$.
Every elliptic subset $\E\subset P(\M)$ defines a finite face
$\gamma (\E) = (\cap_{\delta \in \E}{\Ha_\delta})\cap \M$ of $\M$ of
the codimension $\# \E$. Any finite face $\gamma$ of $\M$ defines an
elliptic subset
$$
P(\gamma^\perp ) =\{\delta \in P(\M)\mid \gamma \subset \Ha_\delta \}
$$
with $\# \E = \codim \gamma$ elements.

A subset $\Pa \subset P(\M)$ is called {\it connected parabolic} if
$\Pa$ has negative  semi-definite (i.e. with no positive
and at least one zero square)
Gram matrix and the Gram graph of $\Pa$ is connected.
One can show that the rank of the Gram matrix of $\Pa$ is equal to
$\# \Pa -1$, and it has exactly one zero square. Connected parabolic
subsets of $P(\M)$ correspond to some infinite points of $\M$. If
$\Pa$ is a connected parabolic subset, then
$\cap_{\delta \in \Pa}{\Ha_\delta}$ is an infinite point of $\M$.
By this construction, two different connected parabolic subsets
$\Pa_1, \Pa_2 \subset P(\M)$ define the same infinite point of
$\M$ if and only if subgraphs $\Gamma (\Pa_1)$ and $\Gamma(\Pa_2)$
are disjoined in $\Gamma (P(\M))$ (i.e. $\Pa_1 \cdot \Pa_2 =0$).
If for a subset $\Pa^\prime\subset P(\M)$ the intersection
$\cap_{\delta \in \Pa^\prime}{\Ha_\delta}$ is an infinite point, then
this point belongs to $\M$ and one of connected components of Gram
graph of $\Pa^\prime$ is parabolic.

A subset $T\subset P(\M)$ is called {\it hyperbolic} if its Gram
matrix has at least one positive square (then it has exactly one
positive square). A subset $L \subset P(\M)$ is called
{\it Lanner} if it is minimal hyperbolic.
Thus, Gram matrix of $L$ is hyperbolic,
but Gram matrix of any proper subset of $L$ is not
hyperbolic (i.e. it is negative definite or negative semi-definite).
It is not difficult to see that
Gram graph $\Gamma (L)$ is connected (see Lemma 1.5.1.2 below).
Any proper subset of $L$ is
either elliptic or connected parabolic.

Obviously, a subset
$T\subset P(\M)$ is hyperbolic
if and only if $T$ contains a Lanner subset.
It is clear that a subset $T\subset P(\M)$ is not hyperbolic
if and only if any
connected component of $\Gamma (T)$ is either elliptic or
connected parabolic. In particular, by properties of elliptic and
connected parabolic subsets mentioned above, a non-hyperbolic
subset $T$ is orthogonal to some point of $\M$. The last property is
very important. We have

\proclaim{Lemma 1.3.1} A polyhedron $\M$ with acute angles
in a Lobachevsky space $\La$
is finite in infinite points of $\La$.
\endproclaim

\demo{Proof} Let $P=\r^+c \in \La$
where $c\in \Phi$ and $c\not=0$ but
$c^2=0$. The set
$$
Q = \{ \delta \in P(\M)\mid \delta \cdot c=0 \}
$$
is a subset of $c^\perp$. Since $\Phi$ is a hyperbolic space, the
linear subspace $c^\perp \subset \Phi$
is negative semi-definite. It follows
that Gram matrix of $Q$ is not hyperbolic. Thus, connected
components of $\Gamma (Q)$ are either negative definite or
negative semi-definite. By results mentioned above,
$\# Q \le 2\rk (\Gamma (Q)) \le 2(\dim \La +1)$. In particular,
$Q$ is finite.
$\blacksquare$
\enddemo

\subhead
1.4. Diagram method
\endsubhead

Series of articles of the author \cite{10}, \'E.B. Vinberg
\cite{6},
A.G. Khovanskii \cite{15} and M.N. Prokhorov \cite{14} was
devoted to getting a bound for dimension of Lobachevsky space $\La$
admitting a reflection group (i.e. a discrete
group generated by reflections
in hyperplanes of $\La$) with a fundamental polyhedron $\M$
of a finite volume. It was shown that $\dim \La <996$.
We remarked in \cite{11}, \cite{12} that in fact, in these papers,
there was obtained a general result valid for arbitrary
elliptic polyhedra with acute angles in a Lobachevsky space, and
applied these results to Del Pezzo surfaces with log terminal
singularities.

Here we want to extend these results to parabolic and hyperbolic
polyhedra with acute angles. To apply results of Section 1.2, we prove
a general diagram method estimate
for dimension of an elliptic face $\gamma$ of
a polyhedron $\M$ with acute angles.

\proclaim{Theorem 1.4.0} Let $\M$ be a non-degenerate locally-finite convex
polyhedron with acute angles in a Lobachevsky space $\La$,
and $\gamma$ an elliptic face of $\M$.
Let
$$
P(\gamma, \M)=\{\delta \in P(\M) \mid
\overline \Ha_\delta \cap \overline\gamma \not=\emptyset\}
$$
and
$$
P(\gamma^\perp )=\{\delta \in P(\M) \mid \gamma \subset \Ha_\delta \}.
$$
Assume
that there are some
constants $d, C_1, C_2$ such that the conditions (a) and (b)
below hold:

(a)
For any Lanner subset
$L \subset P(\gamma , \M)$ such that $L$ contains at
least two elements
which don't belong to $P(\gamma^\perp )$ and for any proper subset
$L^\prime \subset L$ the set $P(\gamma^\perp)\cup L^\prime$ is
not hyperbolic,
$$
\text{diam}~\Gamma (L)\le d.
$$

(b) For any elliptic subset $\E$ such that
$P(\gamma^\perp )\subset \E \subset P(\gamma )$ and $\E$ has $\dim \La -1$
elements, we have for the distance in
the graph $\Gamma (\E)$:
$$
\sharp \{ \{ \delta_1, \delta_2 \}\subset \E-P(\gamma^\perp)
\mid 1 \le \rho (\delta_1,\delta_2)\le d\}
\le C_1\sharp (\E-P(\gamma^\perp ));
$$
and
$$
\sharp \{ \{\delta_1, \delta_2\} \subset \E-P( \gamma^\perp )
\mid d+1\le \rho (\delta_1,\delta_2) \le 2d+1\}
\le C_2 \sharp (\E-P(\gamma^\perp )).
$$
Then $\dim ~\gamma <96(C_1+C_2/3)+68$.

\endproclaim

\demo{Proof} See the proof in Sect. 1.5.
$\blacksquare$
\enddemo

In particular, for elliptic $\M$, we get the following statement which is
a variant of the similar statement we had formulated in
\cite{11} and \cite{12} (in \cite{11} and \cite{12} we used
$\sharp L$ instead of $\text{diam}~ \Gamma (L)$).

\proclaim{Theorem 1.4.1} Let $\M$ be an elliptic convex
polyhedron with acute angles in a Lobache\-vsky space $\La$.
Assume that there are some constants $d, C_1, C_2$
such that the conditions (a) and (b) below hold:

(a)
For any Lanner subset
$L \subset P(\M)$
$$
\text{diam}~\Gamma (L)\le d.
$$

(b) For any elliptic subset $\E \subset P(\M)$
such that $\E$ has $\dim \La -1$
elements, we have for the distance in
the graph $\Gamma (\E)$:
$$
\sharp \{ \{ \delta_1, \delta_2 \}\subset \E
\mid 1 \le \rho (\delta_1,\delta_2)\le d\}
\le C_1\sharp \E;
$$
and
$$
\sharp \{ \{\delta_1, \delta_2\} \subset \E
\mid d+1\le \rho (\delta_1,\delta_2) \le 2d+1\}
\le C_2 \sharp \E.
$$
Then $\dim \La <96(C_1+C_2/3)+68$.
\endproclaim

For a parabolic (relative to an infinite point $P\in \La_\infty$) or
hyperbolic (relative to a subspace $\T \subset \La$) polyhedron
$\M \subset \La$ with acute angles,
by Lemmas 1.2.1 and 1.2.2, there exists
an elliptic face $\gamma \subset \M$. We have mentioned in Sect. 1.3 that
$\gamma$ has acute angles (see the proof of Theorems 1.4.2 and 1.4.3
below).
There are two possibilities to apply results above to get an estimate for
$\dim \La$ using this elliptic face $\gamma$.
In the first place, we can apply Theorem 1.4.1 to the elliptic
polyhedron  $\gamma$ with acute angles.
Secondly, we can directly apply Theorem 1.4.0 to the pair
$\gamma \subset \M$.
Both these variants give good results which are almost the same in
practice. We formulate results of both.

At first, we apply Theorem 1.4.1. For a subset $Q$ of
the set of vertices of a graph $\Gamma$ and
vertices $v_1,v_2$ of $\Gamma$ which both don't belong to  $Q$,
we consider the distance
$$
\rho_Q(v_1,v_2) =\min \{\min_s{(\rho (s)-\#(s\cap Q))}, +\infty \},
$$
where $s$ is a path in $\Gamma$ joined $v_1, v_2$, the distance
$\rho (s)$ is the length of $s$ in $\Gamma$, and the $\#(s\cap Q)$ is
the number of vertices of $s$ which belong to $Q$.

\proclaim{Theorem 1.4.2} Let $\M$ be a convex parabolic
relative to a $P \in \La_\infty$
polyhedron with acute angles in a Lobachevsky space $\La$.
Let $\gamma$ be the elliptic codimension one face of
$\M$ from Lemma 1.2.1, and
$\gamma \subset \Ha_e$, $P \notin \overline\Ha_e$ where $e\in P(\M)$.
Assume that there are some constants $d, C_1, C_2$
such that the conditions (a) and (b) below hold:

(a)
For any Lanner subset
$L \subset P(\gamma , \M)$ such that $L$ contains at
least two elements different from $e$ and for any proper subset
$L^\prime \subset L$ the set $\{e\} \cup L^\prime$ is
not hyperbolic,
$$
\text{diam}~\Gamma (L)\le d.
$$

(b) For any elliptic subset $\E \subset P(\M)$
such that $\E$ has $\dim \La -1$
elements and $e\in \E$, we have for the distance in
the graph $\Gamma (\E)$:
$$
\sharp \{ \{ \delta_1, \delta_2 \}\subset \E-\{e\}
\mid 1 \le \rho_e (\delta_1,\delta_2)\le d\}
\le C_1(\sharp \E-1);
$$
and
$$
\sharp \{ \{\delta_1, \delta_2\} \subset \E-\{e\}
\mid d+1\le \rho_e (\delta_1,\delta_2) \le 2d+1\}
\le C_2 (\sharp \E - 1).
$$
Then $\dim \La <96(C_1+C_2/3)+69$.
\endproclaim

\proclaim{Theorem 1.4.3} Let $\M$ be a convex hyperbolic relative
to a subspace $\T \subset \La$
polyhedron with acute angles in a Lobachevsky space $\La$.
Let $\gamma$ be the elliptic and of codimension $\le \dim \T$
face of $\M$ from Lemma 1.2.2, and
$Q=P(\gamma ^\perp)$ where $Q\subset P(\M)$, $Q$ has
$\codim \gamma \le \dim \T$
elements and $\T \not\subset \Ha_\delta$ for any $\delta\in Q$.
Assume that there are some constants $d, C_1, C_2$
such that the conditions (a) and (b) below hold:

(a)
For any Lanner subset
$L \subset P(\gamma ,\M )$ such that $L$ contains at
least two elements which don't belong to $Q$ and for any proper subset
$L^\prime \subset L$ the set $Q \cup L^\prime$ is not hyperbolic,
$$
\text{diam}~\Gamma (L)\le d.
$$

(b) For any elliptic subset $\E\subset P(\M)$
such that $\E$ has $\dim \La -1$
elements and $Q \subset \E$, we have for the distance in
the graph $\Gamma (\E)$:
$$
\sharp \{ \{ \delta_1, \delta_2 \}\subset \E-Q
\mid 1 \le \rho_Q (\delta_1,\delta_2)\le d\}
\le C_1\sharp (\E-Q);
$$
and
$$
\sharp \{ \{\delta_1, \delta_2\} \subset \E-Q
\mid d+1 \le \rho_Q (\delta_1,\delta_2) \le 2d+1\}
\le C_2 \sharp (\E - Q).
$$
Then $\dim \La <96(C_1+C_2/3)+68 + \dim \T $.
\endproclaim

\demo{Proof of Theorems 1.4.2 and 1.4.3}  Let $e\in \Phi$ and
$e^2<0$. Let us consider the orthogonal projection of $\Phi$
along $e$. By this projection, the image $a^\prime$ of $a\in \Phi$ is
equal to
$$
a^\prime=a-e(e\cdot a)/e^2.
\tag4.1
$$
For images $a^\prime, b^\prime$ of $a, b\in \Phi$, we then get
$$
a^\prime \cdot b^\prime =a\cdot b + (e\cdot a)(e\cdot b)/(-e^2).
\tag4.2
$$
We mention, that using these formulae, one can easily deduce that
a face of a convex polyhedron with acute angles is a polyhedron
with acute angles too.

We put $Q=\{e\}$ for the Theorem 1.4.2. Thus, for both
Theorems 1.4.2 and 1.4.3, $Q=P(\gamma^\perp)$ and
$\#Q=\codim \gamma$ (by results mentioned in Sect. 1.3).
The set $P(\gamma)$ of the polyhedron $\gamma \subset [\gamma]$
is the orthogonal projection along $Q$ of the set
$$
P^\prime (\gamma, \M)=\{\delta \in P(\M) \mid Q\cup \{\delta\}\
\text{is an elliptic subset of }P(\M) \}.
$$
In particular, using (4.2) $\# Q$ times, we then get that $\gamma$ is
a polyhedron with acute angles. A subset
$U \subset P(\gamma)$ is elliptic, parabolic or hyperbolic if and
only if it is a projection along $Q$ of an elliptic,
parabolic or hyperbolic subset respectively of
$U^\prime \subset P^\prime (\gamma, \M)$ such that $Q\subset U^\prime$.
Let $a, b\in U^\prime-Q$ and $a^\prime, b^\prime$ their images in
$U$. Applying (4.2) several times,  we get
$$
\rho(a^\prime, b^\prime)=\rho_Q(a, b)
\tag4.3
$$
where $\rho(a^\prime, b^\prime)$ is the distance in $\Gamma (U)$
and $\rho_Q(a,b)$ uses the distance in $\Gamma (U^\prime )$.

Now, we apply Theorem 1.4.1 to the elliptic polyhedron
$\gamma \subset [\gamma]$ and the constants $d, C_1,C_2$ of the
corresponding Theorem 1.4.2 or 1.4.3.
By (4.3), the condition (b) of Theorem 1.4.1 is equivalent
to the condition (b) of the corresponding Theorem 1.4.2 or 1.4.3.

Now, let $R\subset P(\gamma)$ be a Lanner subset. Then $R$ is the
image by the projection
of the hyperbolic subset $M\subset P^\prime(\gamma , \M)$
such that $M$ contains $Q$, $M$ contains at least two
elements which don't belong to $Q$ (because $R$ has at least two
elements), and any proper subset $M^\prime \subset M$
which contains $Q$ is not hyperbolic. Since $M$ is hyperbolic,
there exists a Lanner subset $L \subset M$. Let
$T \subset R$ be the image of $L$ by the projection.
Since $T$ is hyperbolic and $R$ is Lanner, we get $T=R$.
By (4.3),
$$
\text{diam}~\Gamma (R)\le \text{diam}~\Gamma (L).
$$
Thus, the condition (a) of the corresponding Theorem 1.4.2 or
1.4.3 implies the condition (a) of Theorem 1.4.1 for $\gamma$.
Then, we can apply
Theorem 1.4.1 to $\gamma \subset [\gamma]$ with constants
$d, C_1, C_2$ of the corresponding Theorem 1.4.2 or 1.4.3 to
get the estimate for $\dim [\gamma]$.

By Lemma 1.2.1 for the parabolic case,
$\dim \La=\dim [\gamma] + 1$. By Lemma 1.2.2, for the
hyperbolic case, $\codim [\gamma]=\dim \La -\dim [\gamma] \le \dim \T$,
and $\dim \La \le \dim [\gamma] +\dim \T$. Thus, by Theorem 1.4.1, we
get Theorems 1.4.2 and 1.4.3.
 $\blacksquare$

\enddemo

Applying Theorem 1.4.0 and Lemmas 1.2.1 and 1.2.2, we get

\proclaim{Theorem 1.4.2'} Let $\M$ be a convex parabolic
relative to a $P \in \La_\infty$
polyhedron with acute angles in a Lobachevsky space $\La$.
Let $\gamma$ be the codimension one face of $\M$ from Lemma 1.2.1, and
$\gamma \subset \Ha_e$, $P \notin \overline\Ha_e$ where $e\in P(\M)$.
Assume that there are some constants $d, C_1, C_2$
such that the conditions (a) and (b) below hold:

(a)
For any Lanner subset
$L \subset P(\gamma ,\M)$ such that $L$ contains at
least two elements different from $e$ and for any proper subset
$L^\prime \subset L$ the set $\{e\} \cup L^\prime$ is
not hyperbolic,
$$
\text{diam}~\Gamma (L)\le d.
$$

(b) For any elliptic subset $\E\subset P(\M)$
such that $\E$ has $\dim \La -1$
elements and $e\in \E$, we have for the distance in
the graph $\Gamma (\E)$:
$$
\sharp \{ \{ \delta_1, \delta_2 \}\subset \E-\{e\}
\mid 1 \le \rho(\delta_1,\delta_2)\le d\}
\le C_1(\sharp \E-1);
$$
and
$$
\sharp \{ \{\delta_1, \delta_2\} \subset \E-\{e\}
\mid d+1\le \rho (\delta_1,\delta_2) \le 2d+1\}
\le C_2 (\sharp \E - 1).
$$
Then $\dim \La <96(C_1+C_2/3)+69$.
\endproclaim

\proclaim{Theorem 1.4.3'} Let $\M$ be a convex hyperbolic
relative to a subspace $\T \subset \La$
polyhedron with acute angles in a Lobachevsky space $\La$.
Let $\gamma$ be the codimension $\le \dim \T$
face of $\M$ from Lemma 1.2.2, and
$Q=P(\gamma ^\perp)$ where $Q\subset P(\M)$, $Q$ has
$\codim \gamma \le \dim \T$
elements and $\T \not\subset \Ha_\delta$ for any $\delta\in Q$.
Assume that there are some constants $d, C_1, C_2$
such that the conditions (a) and (b) below hold:

(a)
For any Lanner subset
$L \subset P(\gamma , \M )$ such that $L$ contains at
least two elements which don't belong to $Q$ and for any proper subset
$L^\prime \subset L$ the set $Q \cup L^\prime$ is not hyperbolic,
$$
\text{diam}~\Gamma (L)\le d.
$$

(b) For any elliptic subset $\E \subset P(\M)$
such that $\E$ has $\dim \La -1$
elements and $Q \subset \E$, we have for the distance in
the graph $\Gamma (\E)$:
$$
\sharp \{ \{ \delta_1, \delta_2 \}\subset \E-Q
\mid 1 \le \rho(\delta_1,\delta_2)\le d\}
\le C_1\sharp (\E-Q);
$$
and
$$
\sharp \{ \{\delta_1, \delta_2\} \subset \E-Q
\mid d+1 \le \rho (\delta_1,\delta_2) \le 2d+1\}
\le C_2 \sharp (\E - Q).
$$
Then $\dim \La <96(C_1+C_2/3)+68 + \dim \T $.
\endproclaim

\subhead
1.5. The proof of Theorem 1.4.0
\endsubhead

\subsubhead
1.5.1. The compact case
\endsubsubhead

For it would be easier for a reader, we at first prove similar statement
for the case when $\gamma$ does not have an infinite vertex. This statement
is important as itself. One can also apply this to reflection groups and
algebraic surfaces. We don't do this, but a reader can find
appropriate applications himself.  We prove

\proclaim{Theorem 1.5.1.1} Let $\M$ be a non-degenerate locally-finite convex
polyhedron with acute angles in a Lobachevsky space,
and $\gamma$ be an elliptic compact (i.e. without infinite
vertices) face of $\M$.
Let
$$
P(\gamma, \M)=\{\delta \in P(\M) \mid
\overline \Ha_\delta \cap \overline\gamma \not=\emptyset\}
$$
and
$$
P(\gamma^\perp )=\{\delta \in P(\M) \mid \gamma \subset \Ha_\delta \}.
$$
Assume
that there are some
constants $d, C_1, C_2$ such that the conditions (a) and (b)
below hold:

(a)
For any Lanner subset
$L \subset P(\gamma ,\M)$ such that $L$ contains at
least two elements
which don't belong to $P(\gamma^\perp )$ and for any proper subset
$L^\prime \subset L$ the set $P(\gamma^\perp)\cup L^\prime$ is
not hyperbolic,
$$
\text{diam}~\Gamma (L)\le d.
$$

(b) For any elliptic subset $\E$ such that
$P(\gamma^\perp )\subset \E \subset P(\gamma )$ and $\E$ has $\dim \La $
elements, we have for the distance in
the graph $\Gamma (\E)$:
$$
\sharp \{ \{ \delta_1, \delta_2 \}\subset \E-P(\gamma^\perp)
\mid 1 \le \rho (\delta_1,\delta_2)\le 2d+1\}
\le C\sharp (\E-P(\gamma^\perp )).
$$
Then $\dim ~\gamma < 8C + 6$.

\endproclaim

\demo{Proof (Compare with \cite{6} and also \cite{34})}

\proclaim{Lemma 1.5.1.2} Let $\M$ be a polyhedron with acute angles.
Then any Lanner subset $L \subset P(\M)$
has a connected graph $\Gamma (L)$, and any two Lanner
subsets $L, N\subset P(\M)$ either have a common element or
are joined by an edge in $\Gamma(L \cup \N)$.
\endproclaim

\demo{Proof} Let $L=L_1\cup L_2$ where $L_1\cdot L_2=0$ (equivalently,
$\Gamma (L)$ is the disjoint union of $\Gamma (L_1)$ and $\Gamma (L_2)$.
Since $L$ is hyperbolic, one of the subsets $L_1, L_2$ is hyperbolic.
If this is $L_1$, one gets $L=L_1$ since $L$ is minimal hyperbolic.

Let $L\cdot N=0$. Since $L$ and $N$ are both hyperbolic, the
$L \cup N$ has Gram matrix with at least two positive squares. We
get a contradiction since Gram matrix of $P(\M)$ has not more than one
positive square.
$\blacksquare$
\enddemo

A convex elliptic (i.e. closed) polyhedron is called
{\it simplicial} if all its proper faces are simplexes. A  convex
elliptic polyhedron is called {\it simple} (equivalently,
it has simplicial angles) if it is dual to a simplicial one.
In other words, any its face of codimension
 $k$ is contained exactly in $k$ faces of
the highest dimension (the last definition is valid not only for
elliptic polyhedron). By results mentioned in Sect. 1.3,
any elliptic compact convex polyhedron with acute angles is simple.
This is the very important combinatorial property of a compact
polyhedron with acute angles.
Thus, the face $\gamma \subset \M$ is a simple polyhedron.

Besides, since any face $\gamma_1\subset \gamma$
is  a finite face (since $\gamma$ has no vertices at infinity)
of the polyhedron $\M$
and $\M$ has acute angles, we also have:
$$
\sharp P(\gamma_1^\perp ) - \sharp P(\gamma ^\perp )=
\text{codim}_\gamma \gamma_1
\tag1.5.1.1
$$
for any face $\gamma_1$ of $\gamma $, and $\sharp P(\gamma_1^\perp )=
\text{codim}_{\M} \gamma_1$.

We use the following Lemma 1.5.1.3 which was proved in \cite{10}. The lemma
was used in \cite{10} to get a bound ($\le 10$) for the dimension
of a Lobachevsky space admitting an
action of an arithmetic reflection group with a field
of definition of the degree $>N$. Here $N$ is some absolute constant.

\proclaim{Lemma 1.5.1.3}
Let ${\Cal N}$ be a convex elliptic simple
polyhedron of a dimension $n$, and $A_n^{i,k}$ the average number of
$i$-dimensional faces of $k$-dimensional faces of $\Cal N$.
Then for $n\ge 2k-1$

$$
A_n^{i,k}< {{n-i\choose n-k}\cdot
({[n/2]\choose i}+{n-[n/2]\choose i})\over
{[n/2]\choose k}+{n-[n/2]\choose k}}.
$$
In particular:

If $n\ge 3$
$$A_n^{0,2} <
\cases
4(n-1)\over n-2  &\text{if $n$ is even,}\\
4n\over n-1 &\text{if $n$ is odd.}
\endcases
$$

If $n\ge 5$
$$
A_n^{1,3} <
\cases
12(n-1)\over n-4  &\text{if $n$ is even,}\\
12n\over n-3 &\text{if $n$ is odd.}
\endcases
$$
and
$$
A_n^{2,3} <
\cases
6(n-2)\over n-4  &\text{if $n$ is even,}\\
6(n-1)\over n-3 &\text{if $n$ is odd.}
\endcases
$$
If $n\ge 7$,
$$
A_n^{3,4} <
\cases
8(n-3)\over n-6  &\text{if $n$ is even,}\\
8(n-2)\over n-5  &\text{if $n$ is odd.}
\endcases
$$

\endproclaim

\demo{Proof} See \cite{10}. We mention that the right side of the
inequality of the Lemma 1.5.1.3 decreases and tends to the number
$2^{k-i}{k\choose i}$ of $i$-dimensional faces of $k$-dimensional
cube if $n$ increases.
$\blacksquare$
\enddemo

 From the estimate of $A_n^{0,2}$ of the Lemma,
it follows the following Vinberg's Lemma from \cite{6}.
This Lemma was used by \'E.B. Vinberg to get
an estimate ($\dim <30$) for the dimension
of a Lobachevsky space admitting an action of
a discrete reflection group with a compact fundamental polyhedron.

By definition, a {\it 2-angle} of a polyhedron $T$ is an angle of a
2-dimensional face of $T$. Thus, the 2-angle is defined by
a vertex $A$ of $T$, a plane containing $A$ and a 2-dimensional face
$\gamma_2$ of $T$, and two rays with the beginning at $A$ which
contain two corresponding sides of the $\gamma_2$.

\proclaim{Lemma 1.5.1.4} Let $\Cal N$ be a convex elliptic
simple polyhedron of a dimension $n$. Let $C$ and $D$ are some numbers.
Assume that
2-angles of $\Cal N$ are supplied with
weights and the following conditions (1) and (2) hold:

(1) The sum of weights of all 2-angles at any vertex of
$\Cal N$ is not greater than $Cn+D$.

(2) The sum of weights of all angles of any $2$-dimensional
face of $\Cal N$ is at least $5-k$ where $k$ is the number of
vertices of the $2$-dimensional face.

Then
$$
n<8C+5+
\cases
1+8D/n &\text{if $n$ is even,}\\
(8C+8D)/(n-1) &\text{if $n$ is odd}
\endcases .
$$
In particular, for $C\ge 0$ and $D=0$, we have
$$
n<8C+6.
$$
\endproclaim

\demo{Proof}
Since the proof of Lemma is very short, we give the proof here.

Let $\Sigma$ be the sum of weights of all 2-angles of the polyhedron
$\Cal N$. Let $\alpha _0$ be the number of vertices of
$\Cal N$ and $\alpha _2$ the number of $2$-dimensional faces
of $\Cal N$. Since $\Cal N$ is simple,
$$
\alpha _0{n(n-1)\over 2}=\alpha _2 A_n^{0,2}.
$$
 From this equality and conditions of the Lemma, we get inequalities
$$
(Cn+D)\alpha _0\ge \Sigma \ge \sum \alpha _{2,k}(5-k)=
5\alpha _2-\alpha _2A_n^{0,2}=
$$
$$
=\alpha _2(5-A_n^{0,2})=
\alpha_0(n(n-1)/2)(5/A_n^{0,2}-1).
$$
Here $\alpha _{2,k}$ is the number of 2-dimensional faces  with
$k$ vertices of $\Cal N$.  Thus, from this inequality
and the inequality for $A_n^{0,2}$ of Lemma 1.5.1.3, we get
$$
Cn+D\ge (n(n-1)/2)(5/A_n^{0,2}-1)>
\cases
n(n-6)/8 &\text{if $n$ is even,}\\
(n-1)(n-5)/8 &\text{if $n$ is odd.}
\endcases
$$
 From this calculations, Lemma 1.5.1.4 follows.
$\blacksquare$
\enddemo

The proof of Theorem 1.5.1.1. (Compare with \cite{6} and also \cite{34}.)
Let $\angle$ be a 2-angle of $\gamma $.
Let $P(\angle)\subset P(\gamma )$ be
the elliptic set of all $\delta \in P(\M)$ which are orthogonal
to the vertex $\r^+h\in \La$ of the $\angle$.
By (1.5.1.1), the set
$P(\angle)$ is the disjoint union
$$
P(\angle)=P(\angle^\perp)\cup \{ \delta_1(\angle)\} \cup \{ \delta_2(\angle)\}
$$
where $P(\angle^\perp)$ contains all $\delta \in P(\M)$
orthogonal to the plane of
the 2-angle $\angle$, and $\delta_1(\angle )$ and $\delta_2(\angle )$
are orthogonal to two sides of the 2-angle $\angle$.
And by (1.5.1.1), $\sharp P(\angle)=\dim \La$, and elements
$\delta_1(\angle), \delta_2(\angle)$ are defined uniquely.
Backward, an  elliptic subset $P(\angle)\subset P(\M)$ with
$\dim \La$ elements
and a pair of its different elements
$\{\delta_1(\angle ), \delta_2(\angle )\} \subset P(\angle )$
uniquely define the 2-angle $\angle$ of $\M$.
We define the weight $\sigma (\angle )$ by the formula:
$$
\sigma (\angle)=
\cases
1, &\text{if $1\le \rho (\delta_1(\angle ), \delta_2(\angle ))\le 2d+1$,}\\
0, &\text{if $2d+2\le \rho (\delta_1(\angle ), \delta_2(\angle ))$.}
\endcases
$$
Here we take the distance in the graph $\Gamma (P(\angle))$.
Let us prove conditions of the Lemma 1.5.1.4 with the constant
$C$ of Theorem 1.5.1.1 and $D=0$.

The condition (1) follows from the condition (b) of the theorem. We remark
that $\delta_1(\angle), \delta_2(\angle)$ do not belong to the set
$P(\gamma^\perp)$ since $P(\gamma^\perp)\subset P(\angle^\perp)$.

Let us prove the condition (2).

Let $\gamma _3$ be a 2-dimensional triangle face (triangle) of
$\gamma $. The set $P(\gamma _3)$
of all $\delta\in P(\M)$ orthogonal
to some points of $\gamma_3$ is the union of the set
$P(\gamma _3^\perp)$ of $\delta \in P(\M)$, which
are orthogonal to the plane of the triangle
$\gamma _3$,
and $\delta_1, \delta_2, \delta_3$, which are orthogonal to the sides of
the triangle $\gamma _3$. Union of the set
$P(\gamma _3^\perp)$ with any two elements from
 $\delta_1, \delta_2, \delta_3$  is elliptic, since it is orthogonal to
a vertex of $\gamma _3$.
On the other hand, the set $P(\gamma _3)=
P(\gamma _3^\perp)\cup \lbrace \delta_1, \delta_2, \delta_3 \rbrace$
is hyperbolic, since it is not orthogonal to a point of $\M$.
Indeed, the set of all points of $\M$, which are orthogonal
 to the set
$P(\gamma _3^\perp)\cup \lbrace \delta_2, \delta_3 \rbrace$,
$P(\gamma _3^\perp)\cup \lbrace \delta_1, \delta_3 \rbrace$,
or
$P(\gamma _3^\perp)\cup \lbrace \delta_1, \delta_2\rbrace$ is the
vertex $A_1, A_2$ or $A_3$ respectively of the triangle
$\gamma _3$, and
the intersection of these sets of vertices is empty.
Thus, there exists a Lanner subset $L \subset P(\gamma _3)$,
which contains the set $\{ \delta_1, \delta_2, \delta_3\}$.
By the condition (a), the graph
$\Gamma (L)$
contains a shortest path $s$ of the
length $\le d \le 2d+1$ which connects
$\delta_1, \delta_3$.
If this path does not contain $\delta_2$, then the
angle of $\gamma _3$ defined by the set
$P(\gamma _3^\perp )\cup \lbrace \delta_1, \delta_3 \rbrace$ and
the pair $\{ \delta_1, \delta_3\}$ has the weight 1.
If this path contains $\delta_2$, then the angle of
$\gamma _3$ defined by the
set
$P(\gamma _3^\perp )\cup \lbrace \delta_1, \delta_2 \rbrace$ and
the pair $\{\delta_1, \delta_2\}$ has the weight $1$.
Thus, we proved that the side
$A_2A_3$ of the triangle $\gamma _3$ defines an angle of
the triangle with the weight $1$ and one
side $A_2A_3$ of the angle. The triangle has three sides. Thus,
there are at least two angles of the triangle with the
weight $1$.
It follows the condition (2) of the Lemma 1.5.1.4 for the triangle.

Let $\gamma _4$ be a 2-dimensional quadrangle face (quadrangle)
of $\gamma$. In this case,
$$
P(\gamma _4)=
P(\gamma _4^\perp)\cup \lbrace \delta_1, \delta_2, \delta_3, \delta_4
\rbrace
$$
where $P(\gamma _4^\perp)$ is the set of all $\delta \in P(\M)$
which are orthogonal to the plane of the quadrangle
and $\delta_1, \delta_2, \delta_3, \delta_4$ are orthogonal to
the consecutive sides of the quadrangle. As above, one can see that the
sets
$P(\gamma _4^\perp )\cup \lbrace \delta_1, \delta_3\rbrace$,
$P(\gamma _4^\perp )\cup \lbrace \delta_2, \delta_4\rbrace$
are hyperbolic, but the sets
$P(\gamma _4^\perp )\cup \lbrace \delta_1,\delta_2\rbrace$,
$P(\gamma _4^\perp )\cup \lbrace \delta_2,\delta_3\rbrace$,
$P(\gamma _4^\perp )\cup \lbrace \delta_3,\delta_4\rbrace$,
and
$P(\gamma _4^\perp ) \cup \lbrace \delta_4,\delta_1\rbrace$
are elliptic. It follows that there are Lanner subsets $L, N$
such that
$\lbrace \delta_1, \delta_3 \rbrace \subset L \subset
P(\gamma _4^\perp )\cup \lbrace \delta_1, \delta_3\rbrace$
and
$\lbrace \delta_2, \delta_4 \rbrace \subset N \subset
P(\gamma _4^\perp)\cup \lbrace \delta_2, \delta_4\rbrace$.
By Lemma 1.5.1.2, there
exist $\alpha \in L$ and $\beta \in N$ such  that $\alpha\beta$ is
an edge in $\Gamma (\M)$. By the condition (a) of the theorem,
one of $\delta_1, \delta_3$ is
joined by a path $s_1$ of the length $\le d$ with $\alpha$
and this path does not contain another element from $\delta_1,\delta_3$
(here $\alpha$ is the terminal of the path $s_1$).
We can assume  that this
element is $\delta_1$ (otherwise, one
should replace the $\delta_1$ by the $\delta_3$).  As above,
we can assume that the element $\beta$ is
joined by the path $s_2$ of the length $\le d$ with the
element  $\delta_2$ and this path does not contain the element
$\delta_4$. The path
 $s_1 \alpha \beta s_2$ is a path of the length $\le 2d+1$ in the graph
$\Gamma (P(\gamma _4^\perp))\cup \lbrace \delta_1,\delta_2\rbrace)$.
It follows that the angle of the
quadrangle $\gamma _4$, such that two
sides of this angle are orthogonal to the elements
$\delta_1$ and $\delta_2$,
has the weight $1$. It proves the condition (2)
of the Lemma 1.5.1.4 and the theorem.
$\blacksquare$
\enddemo

\subsubhead
1.5.2. The non-compact case
\endsubsubhead

\demo{The proof of Theorem 1.4.0} This proof
is very similar to the proof of Theorem 1.5.1.1
above but details are more complicated. And we
only outline the proof. One can find omitted details
in A.G Khovanskii \cite{15} and M.N. Prokhorov
\cite{14}: we don't want to rewrite these papers here.

In the first place, for the general case when $\gamma$ is only
an elliptic polyhedron with acute angles with some vertices at
infinity, the $\gamma$ may not be a simple polyhedron. But it is
"almost" simple.

A convex elliptic (i.e. closed) polyhedron is called
{\it simplicial in codimension $k$}
if all its proper faces of codimension $\ge k+1$ are simplexes.
Thus, a simplicial polyhedron is simplicial in codimension $0$.
A  convex closed polyhedron is called {\it simple in dimension
$k$} (equivalently, it has simplicial angles in dimension $k$)
if it is dual to a polyhedron which is simplicial in codimension $k$.
In other words, a polyhedron $T$ is simple in dimension $k$ if
any its proper face of a dimension $m \ge k$
is contained exactly in $\dim T - m$ faces of
highest dimension (the last definition is valid not only
for an elliptic polyhedron). By results mentioned in Sect. 1.3,
any convex polyhedron with acute angles is simple in dimension
$1$. This is the very important combinatorial property of a
polyhedron with acute angles.
Thus, the face $\gamma \subset \M$ is an elliptic simple in
dimension $1$ polyhedron.

Besides, since any face $\gamma_1\subset \gamma$ of
$\dim \gamma_1 \ge 1$ is a finite face of the polyhedron
$\M$ with acute angles, we also have:
$$
\sharp P(\gamma_1^\perp ) - \sharp P(\gamma ^\perp )=
\text{codim}_\gamma \gamma_1
\tag1.5.2.1
$$
for any face $\gamma_1$ of $\gamma $ if $\dim \gamma_1\ge 1$,
and $\sharp P(\gamma_1^\perp)=\text{codim}_{\M}\gamma_1$.

A.G. Khovanskii \cite{15} generalized Lemma 1.5.1.3 to polyhedra
which are simple in dimension $1$.

\proclaim{Lemma 1.5.2.3} Lemma 1.5.1.3 is valid for
convex elliptic polyhedra which are simple in dimension $1$.
\endproclaim

\demo{Proof} See \cite{15}. We should say that the proof of
this Lemma required new ideas comparing with the proof of
Lemma 1.5.1.3 in \cite{10}. We remark that for a rational
polyhedron (actually, only this case is important for surfaces)
one can prove Lemma 1.5.2.3 similarly to \cite{10}
using known properties of the combinatorial polynomial of the
polyhedron.
$\blacksquare$
\enddemo

Since the property (1.5.2.1) is only true for $\dim \gamma_1\ge 1$,
one cannot use $2$-dimensional angles. One has to consider
$3$-dimensional angles (3-angles).

By definition, a {\it 3-angle} of a polyhedron $T$ is a
facial angle of a $3$-dimensional face of $T$.
Thus, the 3-angle is defined by
an edge (i.e. $1$-dimensional face) $\gamma_1$
of $T$, a $3$-dimensional face $\gamma_3$ of
$T$ containing $\gamma_1$ and two different 2-dimensional faces
$\gamma_2^{(1)}, \gamma_2^{(2)}$ of $\gamma_3$
with the common edge $\gamma_1$.

To formulate the analog of Lemma 1.5.1.4, we need
an additional definition.
Let us take two $3$-dimensional tetrahedra and glue together along
their two-dimensional triangle face. This polyhedron is
called triangle bipyramid. Thus, a triangle bipyramid has $5$ vertices:
two vertices of the order $3$
and three vertices of the order $4$ (here the order of a vertex is
equal to the number of edges which contain this vertex).
A $3$-dimensional polyhedron which
is combinatorially equivalent to a triangle bipyramid is called {\it bad}.
A $3$-dimensional polyhedron is {\it good} if it is not a bad one.

Using estimates for
$A_n^{3,4}, A_n^{1,3}, A_n^{2,3}$ of Lemma 1.5.2.3,
M.N.Prokhorov \cite{14}
proved the following analog of Vinberg's Lemma 1.5.1.4.

\proclaim{Lemma 1.5.2.4} Let $\Cal N$ be an elliptic
convex polyhedron with acute angles in a Lobache\-vsky space of
dimension $n$. Let $C>0$ be a constant.
Assume that
3-angles of $\Cal N$ are supplied with
weights and the following conditions (1) and (2) hold:

(1) The sum of weights of all 3-angles at any edge $\gamma_1$ of
$\Cal N$ is not greater than $C(n-1)$.

(2) The sum $\sigma (\gamma_3)$
of weights of all 3-angles of any good $3$-dimensional
face $\gamma_3$ of $\Cal N$ is at least $7-k$ where $k$ is the number of
$2$-dimensional faces of the $3$-dimensional face $\gamma_3$.

Then
$$
n<96C+68.
$$
\endproclaim

\demo{Proof} The proof uses estimates for
$A_n^{3,4}, A_n^{1,3}, A_n^{2,3}$ of Lemma 1.5.2.3 and is
similar to the
proof of Lemma 1.5.1.4. But the proof requires more
calculations.
Mainly, one should analyze how many bad $3$-dimensional
faces may have a $4$-dimensional convex elliptic
polyhedron with acute angles in Lobachevsky space which has $\le 9$ three-
dimensional faces. See details in\cite{14} in
the proof of \cite{14, Lemma 2}.
$\blacksquare$
\enddemo

Now we are ready to discuss the proof of Theorem 1.4.0.

\demo{The proof of Theorem 1.4.0} (Compare with
\cite{14, \S 3}).)
Let $\angle$ be a 3-angle of $\gamma$, and $\angle$ is defined by
faces $\gamma_1\subset \gamma_2^{(1)},\gamma_2^{(2)}\subset \gamma_3$
of $\gamma$ (see the definition above of a 3-angle).
Let $P(\angle)\subset P(\gamma )$ be
the elliptic set of all $\delta \in P(\M)$ which are orthogonal
to the edge $[\gamma_1]$ of the $\angle$.
By (1.5.2.1), the set
$P(\angle)$ is the disjoint union
$$
P(\angle)=P(\angle^\perp)\cup \{ \delta_1(\angle)\} \cup \{ \delta_2(\angle)\}
$$
where $P(\angle^\perp)$ contains all $\delta \in P(\M)$
orthogonal to the $3$-dimensional subspace $[\gamma_3]$ of
the 3-angle $\angle$, and $\delta_1(\angle )$ and $\delta_2(\angle )$
are orthogonal to two sides $[\gamma_2^{(1)}]$ and
$[\gamma_2^{(2)}]$ of the $\angle$.
And, by (1.5.2.1), $\sharp P(\angle)=\dim \La-1$, and elements
$\delta_1(\angle), \delta_2(\angle)$ are defined uniquely.
Backward, an elliptic subset $P(\angle)\subset P(\M)$
with $\dim \La -1$ elements and a pair of its different elements
$\{\delta_1(\angle ), \delta_2(\angle )\} \subset P(\angle )$
uniquely define the 3-angle $\angle$.
We define the weight $\sigma (\angle )$ by the formula:
$$
\sigma (\angle)=
\cases
1, &\text{if $1\le \rho (\delta_1(\angle ), \delta_2(\angle ))\le d$,}\\
1/3, &\text{if $d+1\le \rho (\delta_1(\angle ), \delta_2(\angle ))\le 2d+1$,}\\
0, &\text{if $2d+2\le \rho (\delta_1(\angle ), \delta_2(\angle ))$.}
\endcases
$$
Here we take the distance in the graph $\Gamma (P(\angle))$.
Let us prove conditions of the Lemma 1.5.2.4 with the constant
$C=C_1+C_2/3$ where $C_1,C_2$ are constants from the condition of
Theorem 1.4.0.

The condition (1) follows from the condition (b) of the theorem. We remark
that $\delta_1(\angle), \delta_2(\angle)$ do not belong to the set
$P(\gamma^\perp)$ since $P(\gamma^\perp)\subset P(\angle^\perp)$.

Let us prove the condition (2). Let $\gamma_3\subset \gamma$ be
a good $3$-dimensional face of $\gamma$ and $\gamma_3$ contains
$\le 6$ two-dimensional faces. (We want to show that this case is
similar to the case of triangle or quadrangle for the proof of
Theorem 1.5.1.1.)

A set $\{\gamma^{(1)},...,\gamma^{(k)}\}$ of two-dimensional faces of a
$3$-dimensional face $\gamma_3$
is called {\it good} if it has one of the types 1---4 below
(see \cite{14, Figure 2}):

{\it Type 1}: $k=4$ and
$\gamma^{(1)},...,\gamma^{(4)}$ define the configuration of
two-dimensional faces of a three-dimensional tetrahedron
(then $\gamma_3$ is also tetrahedron).

{\it Type 2}: $k=3$ and
$\gamma^{(1)}\cap \gamma^{(2)}\cap \gamma^{(3)}$ is empty but
$\gamma^{(i)} \cap \gamma^{(j)}$ is an edge of $\gamma_3$ for
any  $1\le i <j \le 3$.

{\it Type 3}: $k=3$ and $\gamma^{(1)}\cap \gamma^{(2)}\cap \gamma^{(3)}$
is empty,
$\gamma^{(1)} \cap \gamma^{(2)}$ and $\gamma^{(2)} \cap \gamma^{(3)}$
are edges of $\gamma_3$ and $\gamma^{(3)} \cap \gamma^{(1)}$
is a vertex at infinity of $\gamma_3$.

{\it Type 4}: $k=4$, intersection of any three different faces
from $\gamma^{(1)},...,\gamma^{(4)}$ is empty,
intersection of any two different faces
from $\gamma^{(1)},...,\gamma^{(4)}$ is empty except
intersections
$\gamma^{(1)} \cap \gamma^{(2)}, \ \gamma^{(2)} \cap \gamma^{(3)},\
\gamma^{(3)} \cap \gamma^{(4)},\
\gamma^{(4)} \cap \gamma^{(1)}$
which define edges of $\gamma_3$.

For a good subset
${\Cal K}=\{\gamma^{(1)},...,\gamma^{(k)}\}$ of two-dimensional faces of
a $3$-dimensional face $\gamma_3$ we consider the sum
$\sigma ({\Cal K})$ of all 3-angles with both sides
which belong to ${\Cal K}$.
Using Lemma 1.5.1.2, like for
the proof of Theorem 1.5.1.1, one can prove that
$$
\sigma ({\Cal K}) \ge
\cases
3, &\text{if ${\Cal K}$ has the type 1},\\
2, &\text{if ${\Cal K}$ has the type 2},\\
1, &\text{if ${\Cal K}$ has the type 3},\\
1/3, &\text{if ${\Cal K}$ has the type 4}.\\
\endcases
\tag1.5.2.2
$$
Here the types 1, 2 and 3 are similar to triangle
and the type 4 is similar to quadrangle.

It is shown in \cite{14} that any good face $\gamma_3$ with
$k\le 6$ two-dimensional faces contains good subsets of two-dimensional
faces. Using these subsets and
(1.5.2.2), one can prove that $\sigma (\gamma_3)\ge 7-k$.
In fact, one should draw pictures of all $3$-dimensional polyhedra
with acute angles and $k \le 6$ (see \cite{14, Figure 3}) to find good
subsets of two-dimensional faces. For example, if $\gamma_3$ is a cube,
$\gamma_3$ evidently has three different good subsets of the type 4.
By (1.5.2.2),
$\sigma (\gamma_3)\ge 3\cdot (1/3)=1=7-k$.
We remark that a bad 3-dimensional face does not have
a good subset of two-dimensional faces.
It is why we had to exclude bad faces in Lemma 1.5.2.4.

This finishes the proof of Theorem 1.4.0.
$\blacksquare$
\enddemo

\enddemo

\head
\S 2. Discrete reflection groups in Lobachevsky spaces and the diagram method
\endhead

Let $G$ be a discrete group in a Lobachevsky space $\La$ with
a fundamental domain of finite volume (i.e. crystallographic group).
Let $W \triangleleft G$ be the reflection subgroup of $G$,
(i.e. $W$ is generated by all reflections in hyperplanes
which belong to $G$).
We choose a fundamental polyhedron $\M$ of $W$. Then
one can consider the quotient group $A=G/W$ as an automorphism
group of $\M$. Thus, $G=W\rtimes A$ is a semi-direct product.
As in \cite{26, Sect. 3}, we say that $G$ is a
{\it generalized (hyperbolic) crystallographic reflection group}
if there exists a subgroup
$A^\prime \subset A$ of a finite index and a non-trivial subspace
$\La^\prime \subset \La$ (including the case when
$\La^\prime$ is an infinite point of $\La$) such that
$\La^\prime$ is $A^\prime$-invariant. Here "non-trivial"
means that $\La^\prime \not=\emptyset$
and $\La^\prime \not= \La$.

In the first place, we have the following preliminary

\proclaim{Proposition 2.1} Let $G$ be a generalized crystallographic
reflection group in a Lobachev\-sky space $\La$. Then there exist
a subgroup $A^\prime \subset A$ of a finite index and a non-trivial
$A^\prime$-invariant subspace $\T \subset \La$ such that $\T$ is
generated by $\T\cap \M$.
\endproclaim

\demo{Proof} If $A$ is finite, $A^\prime$ is a unite subgroup and
$\T$ is any point of $\M$.

Assume that $A$ is infinite. Let $A^\prime \subset A$
be a finite index subgroup
and let $\La^\prime \subset \La$ be a non-trivial $A^\prime$-invariant
subspace. Let $m\in \M$ be a finite point.
Since $A^\prime$ is discrete and infinite,
there is a sequence $a_1,...,a_n,..$ in $A^\prime$ such that
the sequence $a_1(m),...,a_n(m),..$
tends to a point $P\in \M \cap \La_\infty$.
If $\La^\prime = Q$ is an infinite point, any horosphere $\E_{Q,R}$ with
the center $Q$ is $A^\prime$-invariant. It follows that the
sequence $a_1(m),...,a_n(m),...$ tends to $Q$
and $Q=P$ belongs to $\M$.
If $\La^\prime$ is a finite subspace of $\La$, the distances
$\rho (a_i(m), \La^\prime )$ are equal. It follows that
$P\in \overline \La^\prime$.
Thus, at any case, $P\in \M\cap \overline \La^\prime$ and
$\M\cap \overline \La^\prime \not=\emptyset$. Thus, we can set
$\overline \T=[\M\cap \overline \La^\prime]$.
$\blacksquare$
\enddemo

Using Proposition 2.1, we divide generalized crystallographic
reflection groups $G$ in three types.
The group $G$ is called {\it elliptic} if $G$ satisfies Proposition 2.1
with $\T$ which is a finite point $P\in \M$.
The group $G$ is called
{\it parabolic relative to $P$}
if $G$ satisfies Proposition 2.1 with $\T$ which is
an infinite point $P\in \M\cap \La_\infty$.
The group $G$ is called {\it hyperbolic relative to $\T$}
if $G$ satisfies Proposition 2.1 with $\T$  which is a finite subspace
of $\La$ generated by $\T\cap \M$ and $0 < \dim \T < \dim \La -1$.

This definition is justified by the following

\proclaim{Proposition 2.2} Let $G$ be a generalized crystallographic
reflection group in a Lobachev\-sky space $\La$. Then:

(a) $\M$ is elliptic if $G$ is elliptic;

(b) $\M$ is parabolic relative to $P \in \M\cap \La_\infty$ if
$G$ is parabolic relative to $P$.

(c) $\M$ is hyperbolic relative to a subspace $\T \subset \La$ if $G$ is
hyperbolic relative to $\T$.
\endproclaim

\demo{Proof} Considering the subgroup
$G^\prime =W \rtimes A^\prime \subset G$ of a finite index,
we can assume that the subspace $\T\subset \La$ of Proposition 2.1 is
$A$-invariant, and for the case (a), $A$ is trivial and $G=W$.

We use the following well-known
statement about discrete groups $G$ in Lobachevsky space
with a fundamental domain of finite volume (see \cite{27}): the $G$
has a fundamental domain $D$ which is an elliptic (i.e. closed)
convex polyhedron. Besides, we use the following standard statement
about discrete reflection groups: the fundamental polyhedron
$\M$ of the group $W$ is a convex locally-finite polyhedron with acute
angles and facial angles of the form
$\pi/n, n \in {\Bbb N},\ n\ge 2$, and any codimension one
face $\gamma$ of $\M$ defines a reflection in the hyperplane $[\gamma]$
which belongs to $G$. It follows that if
$D\cap \M$ contains an open non-empty subset, then $D\subset \M$ and
$D$ is a fundamental domain for $A$ acting in $\M$.

For the case (a), $D=\M$ is an elliptic polyhedron.

For the case (b), simple standard considerations using facts above
imply that $D=\M \cap C_{\N}$ where $\N\subset \E_P$ is an elliptic polyhedron
on the horosphere $\E_P$ and $\N$ is a fundamental domain for the
action of $A$ in
$$
\M_P = \bigcap\limits_{\delta \in P(\M),\ P\in \Ha_\delta}
{\E_P \cap \Ha_\delta^+} .
$$
Here, by Lemma 1.3.1, the set
$\{\delta \in P(\M) \mid P\in \Ha_\delta \} $ is finite and
the polyhedron
$\M_{P}\subset \E_P$ is bounded by a finite set
of hyperplanes $\Ha_\delta\cap \E_P$. Since $D$ is an elliptic polyhedron,
it easily follows that $\M$ is parabolic relative to $P$.

For the case (c), similarly,
$D=\M \cap C_{\N}$ where $\N\subset \T$ is an elliptic polyhedron in
$\T$ and $\N$ is a fundamental domain for the action of $A$ in
$\M \cap \T$. Since $D$ is an elliptic polyhedron in $\La$, it
easily follows that $\M$ is hyperbolic relative to $\T$.
$\blacksquare$
\enddemo

Using Theorems 1.4.1--1.4.3 (or 1.4.1, 1.4.2', 1.4.3') and
Proposition 2.2, we get

\proclaim{Theorem 2.3} Let $G$ be a generalized crystallographic
reflection group in a Lobachevsky space $\La$. Then

(a) $\dim \La < 1056$ if $G$ is elliptic;

(b) $\dim \La < 1057$ if $G$ is parabolic;

(c) $\dim \La < 1056 + \dim \T$ if $G$ is hyperbolic relative to
a subspace $\T \subset \La$.
\endproclaim

\demo{Proof} It is well-known that $\text{diam} \Gamma (L) \le 8$ for
any Lanner subset $L\subset P(\M)$ (for example, see \cite{14}).
For any elliptic subset $\E\subset P(\M)$, every connected
component of the graph
$\Gamma (\E)$ has the form ${\Bbb A}_n$, ${\Bbb D}_n$,
${\Bbb E}_6$, ${\Bbb E}_7$ or ${\Bbb E}_8$. It easily follows that
we can take $C_1=8,\ C_2=9$ for Theorems 1.4.1--1.4.3
(and 1.4.2', 1.4.3').
$\blacksquare$
\enddemo

\remark{Remark 2.4} With reference to the assertion (c) of
Theorem 2.3, we don't know if there exists an absolute estimate
of $\dim \La$ which does not depend from $\dim \T$.
\endremark

\remark{Remark 2.5} A.G. Khovanskii \cite{15} and M.N. Prokhorov \cite{14}
have proven that $\dim \La < 996$ for the case (a) (see considerations in
\cite{14}). Similarly, we can improve estimates above in this special
case of reflection groups:
$\dim \La < 997$ for the parabolic case, and $\dim \La < 996+\dim \T$
for the hyperbolic case.

If $\M$ is bounded,
\'E.B. Vinberg \cite{6} had proven that $\dim \La < 30$.
Before, the author had proven \cite{10} that $\dim \La < 10$
if $G$ is an arithmetic reflection group (elliptic) and the degree
of the ground field is sufficiently large.
\endremark

\head
\S 3. Diagram method and algebraic surfaces
\endhead

Let $Y$ be a non-singular projective
algebraic surface over an algebraically closed
field and $\NS (Y)$ be the Neron--Severi lattice (i.e. an
integral symmetric bilinear form) of $Y$. By Hodge index
theorem, the lattice $\NS (Y)$ is hyperbolic
(i.e. it is non-degenerate and has exactly one
positive square) and defines the Lobachevsky
space
$\La (Y)= V^+ (Y)/\r^+$  where $V(Y)=\{x\in \NS (Y)\otimes \r \mid
x^2>0\}$ and $V^+(Y)$ is the half-cone containing
the class of a hyperplane section of $Y$.

We recall that a curve $C$ of $Y$ is called {\it exceptional} if
$C$ is irreducible and $C^2<0$. Let $\Exc(Y)\subset \NS(Y)$ be the
set of classes in $\NS(Y)$ of all exceptional curves of $Y$ and
$$
\M(Y)=\bigcap_{\delta \in \Exc(Y)}{\Ha_\delta^+}
$$
the corresponding convex polyhedron in $\La (Y)$ with
$P(\M (Y))=\Exc (Y)$. The $\M(Y)$ is a polyhedron with acute angles
since $\delta \cdot \delta^\prime \ge 0$ if $\delta , \delta^\prime
\in \Exc (Y)$ and $\delta \not=\delta^\prime$.

One can apply the diagram method --- results of Sect. 1.4--- to
$\M(Y)$ (and $Y$) if $\M(Y)\subset \La(Y)$ is locally finite and
is elliptic, parabolic or hyperbolic. Below we consider an
important case when this is true.

We recall that a formal finite linear combination
$D=\sum_i{a_iC_i}$ of irreducible curves $C_i$ of $Y$
is called
{\it divisor, $\Bbb Q$-divisor} and  {\it $\r$-divisor}
if $a_i\in {\Bbb Z}$, $a_i\in {\Bbb Q}$, and $a_i \in \r$
respectively. An $\r$-divisor $D=\sum_i{a_iC_i}$ is called
{\it effective} if all $a_i\ge 0$. An $\r$-divisor $D$
is called numerically effective (equivalently, $nef$) if
$D\cdot F\ge 0$ for any effective divisor $F$. An $\r$-divisor
$F$ is called {\it pseudo-effective} if $F\cdot D\ge 0$ for any
$nef$ divisor $D$. Here it is sufficient considering $D$
which are irreducible curves $D$ with $D^2\ge 0$. The same
definition is valid for classes of divisors
in $N(Y)=\NS(Y)\otimes \r$.

One can also use Mori (or Kleiman--Mori) cone to give these definitions.
Let $NE(Y)\subset N(Y)$ be a convex cone generated by all classes
of effective divisors. The {\it Mori cone $\nem (Y)$} is the closure of
$NE(Y)$ in Euclidean topology of $N(Y)$ and is equal to the set
of all pseudo-effective divisors.
The dual cone $\nem (Y)^\ast \subset N(Y)$
(with respect to the intersection pairing)
is the set of all $nef$ divisors. Since $N(Y)$ is a hyperbolic space with
respect to the intersection pairing, we have
$$
\M(Y)=\nem (Y)^\ast /\r^+ .
$$

We recall Zariski decomposition (see \cite{31}, \cite{28}).
Let $D$ be a pseudo-effective (i.e.
$D\in \nem (Y)$) $\Bbb Q$-divisor. Then $D$ has the
{\it Zariski decomposition} (in $N(Y)$):
$$
D\equiv P+N,
\tag1
$$
where

(i) $P$ is a numerically effective (i. e. $nef$) $\Bbb Q$-divisor,

(ii) $N=\sum_{i=1}^{n}{\alpha_iF_i}$ where $\alpha_i\in {\Bbb Q}$
and $\alpha_i>0$, and $F_1,...,F_n$ are exceptional curves with
a negative definite Gram matrix $(F_i\cdot F_j)$, $1\le i,j \le n$;

(iii) $P\cdot F_i=0$ for any $F_i,\ 1\le i \le n$.

It is known that the Zariski decomposition is unique: the $P\in N(Y)$,
the divisor $N$ are unique if one has properties (i), (ii), (iii).
If $D$ is an effective $\Bbb Q$-divisor, Zariski decomposition is
defined for $\Bbb Q$-divisors: $D=P+N$.

Zariski decomposition defines {\it the numerical Kodaira dimension}
$\nu (D, Y)$ of $D$:
$$
\nu (D,Y)=
\cases
2,\  &\text{if $P^2>0$},\\
1,\  &\text{if $P^2=0$ and $P\not\equiv 0$},\\
0,\ &\text{if $P\equiv 0$}\\
-\infty \ &\text{if $D$ is not pseudo-effective}.
\endcases
$$
Let $K=K_Y$ be the canonical class of $Y$. Then $\nu (K,Y)$ is
called  the {\it numerical Kodaira dimension of $Y$} and
$\nu (-K, Y)$ is called
the {\it numerical anti-Kodaira dimension of $Y$}.

We use Mori theory \cite{23} to prove

\proclaim{Proposition 3.1} Let $Y$ be a non-singular projective
algebraic surface
over an algebraically closed field. Let $K=K_Y \not\equiv 0$,
$\nu (-K, Y)\ge 0$ and
$$
-K\equiv P+N
$$
be the Zariski decomposition of $-K$ where
$N=\sum_{i=1}^{n}{\alpha_iF_i}$. Then we have statements
(a), (b), (c) and (d) below:

(a) All exceptional curves of $Y$ are curves $F_1,...,F_n$, non-singular
rational curves with square $(-2)$ and exceptional curves of the first
kind (i.e. non-singular rational curves with square $(-1)$).
The polyhedron $\M(Y)$ is locally finite.

(b) $\M(Y)$ is elliptic and $\r^+P\in \M(Y)$ if $\nu (-K,Y)=2$ (i.e.
$P^2>0$).

(c) $\M(Y)$ is parabolic relative to the point
$\Pa=\r^+P \in M(Y)\cap \La_\infty$
if $\nu (-K,Y)=1$ (i.e. $P^2=0$ but $P\not\equiv 0$).

(d) $\M(Y)$ is hyperbolic relative to the subspace
$\T=\bigcap_{i=1}^{i=n}{\Ha_{F_i}}$ of $codim \T=n$ in
$\La (Y)$ if $\nu (-K,Y)=0$ (i.e. $P\equiv 0$) .

\endproclaim

\demo{Proof} (a) Let $L$ be an exceptional curve of $Y$ and
$L$ is different from $F_1,...,F_n$. The arithmetic genus
$p_a(L)=(L^2+L\cdot K)/2+1\ge 0$. Since $L$ is different from
$F_1,...,F_n$, then $L\cdot K=-L\cdot P -L\cdot N\le 0$. Hence,
either $L^2=-2$ and $L\cdot K=0$ or
$L^2=L\cdot K=-1$. For both cases $p_a(L)=0$.
Then $p_g(L)=0$ since $p_a(L)\ge p_g(L)\ge 0$, and $L$ is non-singular
rational. This proves the first statement of (a).

Elements $\delta \in \NS (Y)$ with $\delta ^2=-1, -2$ define
reflections
$$
s_\delta: x \mapsto x-2(x\cdot \delta )/\delta^2,\ \ x \in \NS(Y).
$$
Evidently, $s_\delta \in O_+(\NS(Y))$ and $s_\delta$ in $\La(Y)$ is
the reflection in the hyperplane $\Ha_\delta$.
The group $O_+(\NS(Y))$ is discrete in $\La (Y)$ (this is true for
the automorphism group of any hyperbolic lattice $S$ acting in
$\La(S)$, for example, see \cite{27} for the general statement).
It follows that the set of hyperplanes $\Ha_\delta,\ \delta \in \NS(Y)$
with $\delta^2=-1$, or $-2$, is locally finite in $\La(Y)$. The set of the
hyperplanes $\Ha_{F_1},..., \Ha_{F_n}$ is finite. It follows that
$\M(Y)$ is locally finite.

Let us prove (b), (c) and (d).
The polyhedron $\M(Y)$ is obviously elliptic if
$\rho (Y)=\dim \NS(Y)=1$ or $2$ because $\La(Y)$ is $0$ or $1$-dimensional
respectively.

Assume $\rho (Y)\ge 3$. Let
$$
\nem _-(Y) = \{ x\in \nem (Y) \mid x\cdot K<0\}
$$
be the "negative" part of Mori cone $\nem (Y)$.
By Mori theory \cite{23},  the projectivization
$P\nem_-(Y)=\nem _-(Y)/\r^+$ is locally polyhedral and has as its
vertices (equivalently, extremal rays of the cone
$\nem_-(Y)$) rays $\r^+E$ where $E$ is an exceptional curve of the first
kind.

By Riemann--Roch Theorem for surfaces, we have
$$
\overline{V^+(Y)} \subset \nem (Y).
$$
Let $P\nem (Y)= (\nem (Y)-\{0\})/\r^+$ is the projectivization of the
Mori cone. Obviously, the cone $\overline{V^+(Y)}$ is
self-dual: $\overline{V^+(Y)}^{\ \ast} =\overline{V^+(Y)}$.
Thus, we get the sequence of embeddings of projectivizations of
the convex cones:
$$
\M(Y)\subset \overline{\La (Y)}\subset P\nem (Y)
\tag2
$$
which is self-dual:
$$
\M(Y)^\ast =P\nem (Y),\ P\nem(Y)^\ast=\M(Y)\ \text{and}\
\overline{\La(Y)}^{\ \ast}=\overline{\La(Y)}.
\tag3
$$
Using these facts, we get

\proclaim{Lemma 3.2} With conditions of Proposition 3.1,
assume that $\dim \NS(Y)\ge 3$ and
$Q\in \La(Y)_\infty \cap \M(Y)\cap P\nem_-(Y)$.
Let $U$ be an open neighborhood of $Q$ in the infinite
sphere $\La(Y)_\infty$. Then $U$ is not contained in $\M(Y)$.
\endproclaim

\demo{Proof} Let $U\subset \M(Y)$. Considering smaller
neighborhood, we can assume that $U\subset P\nem_-(Y)$
since $\La(Y)_\infty \subset P\nem (Y)$.
We claim that $U$ is contained in the boundary of $P\nem_-(Y)$.
Indeed, let $\r^+c \in U$ where $c\in \NS (Y)\otimes \r$,
$c^2=0$ and $c\not=0$. Since $P\nem (Y)=\M (Y)^\ast$ and
$\r^+c \in \M(Y)$, we have
$$
P\nem_-(Y) \subset P\nem (Y)\subset \overline\Ha_c^+
$$
where
$$
\overline \Ha_c^+=\{ \r^+ x \mid x \in \NS (Y)\otimes \r,\
x \not=0\ \text{and}\ x\cdot c\ge 0\}
$$
is the projectivization of a half-space in $\NS(Y)\otimes \r$. The
point $\r^+c$ belongs to the boundary of $\overline \Ha_c^+$
because $c\cdot c=c^2=0$. It follows that $\r^+c$ belongs to the
boundary of $P\nem (Y)$ and $P\nem_-(Y)$ because
$P\nem_-(Y)\subset P\nem (Y)\subset \overline \Ha_c^+$.
Thus, all points of $U$
belong to the boundary of $P\nem_-(Y)$. The sphere $\La (Y)_\infty$
is strongly convex in $P(\NS(Y)\otimes \r)$. It follows that
any point of $U$ is a vertex of
$P\nem_-(Y)$ (equivalently, an extremal ray of $\nem _-(Y)$.
We get the contradiction with Mori theory.
Thus, $U$ is not contained in $\M(Y)$.
$\blacksquare$
\enddemo

Now we prove (b), (c) and (d). By (a), the set $\Exc (Y)$ is divided on
three subsets:  $\Exc (Y)_3=\{F_1,...,F_n\}$ and $\Exc (Y)_2$,
$\Exc (Y)_1$. The last two contain all
non-singular rational exceptional curves with square $(-2)$ and
$(-1)$ respectively  which do not belong to $\Exc (Y)_3$.
If $F\in \Exc (Y)_2$, then $F\cdot K=F\cdot P =
F\cdot F_i=0$ ($i=1,...,n$). If $E\in \Exc (Y)_1$, we get
$-E\cdot K=E\cdot (P+\sum_{i=1}^{n}{\alpha_iF_i})=1$ where
$E\cdot P\ge 0$ and $E\cdot F_i\ge 0$. In particular,
$\r^+P\in \M(Y)$ for cases (b) and (c).

Let us consider the case (b). Let $\r^+h \in \M(Y)$. Then
$-h\cdot K=h\cdot P+h \cdot \sum_{i=1}^n{\alpha_iF_i}>0$,
because $h\cdot F_i\ge 0$, $\alpha_i>0$ and $h\cdot P>0$ since
$\r^+h\in \overline{\La (Y)}$ and $\r^+P\in \La (Y)$. Thus,
$\M (Y)\subset P\nem _-(Y)$. Moreover, the set
$\Exc (Y)_1\cup \Exc (Y)_2$ is finite. Indeed,
the set of elements $\delta \in \NS(Y)$ with properties
$0\le \delta\cdot P \le 1$ and $\delta^2=-1$ or $-2$ is
finite for a hyperbolic lattice $\NS(Y)$ and the fixed element
$P\in \NS(Y)\otimes {\Bbb Q}$ with $P^2>0$. Thus,
the set $\Exc (Y)$ is finite and $\M(Y)$ is bounded by
a finite set of hyperplanes $\Ha_\delta,\ \delta \in \Exc (Y)$.
If $\M(Y)$ is not elliptic, there exists
a non-empty open subset
$U\subset \M(Y)\cap \La (Y)_\infty$. Since $\M(Y)\subset P\nem_-(Y)$,
we get the contradiction with Lemma 3.2.

Case (c). Similarly to (b), one can prove that
$\M(Y)-\{\Pa\}\subset P\nem_-(Y)$ where $\Pa =\r^+P$. Let
$$
R=\{\delta \in \Exc (Y) \mid \delta\cdot P=0\}.
$$
By Lemma 1.3.1, the set $R$ is finite, since $\M(Y)$ has
acute angles. The set $R$ contains $\{F_1,...,F_n\}$, $\Exc (Y)_2$
and a finite set of elements of $\Exc (Y)_1$. Let
$\Exc (Y)^\prime_1=\Exc (Y)_1-R$. Let
$\M(Y)_\Pa=\bigcap_{\delta \in R} {\overline \Ha_\delta^+}$.
The polyhedron $\M(Y)_\Pa$ is a cone with the vertex $\Pa$ and
the base $T=\{l\in \E_\Pa \mid l\subset \M(Y)_\Pa\}$ on
the horosphere $\E_\Pa$. Obviously, $\M(Y)\subset \M(Y)_\Pa$. Let
${\Cal K}\subset \E_\Pa$ be an elliptic polyhedron on the horosphere.
If $C_{\Cal K}\cap \M(Y)_\Pa =C_{{\Cal K}\cap T}$ is degenerate
(equivalently, the polyhedron ${\Cal K}\cap T\subset \E_\Pa$ is
degenerate), the polyhedron $C_{\Cal K}\cap \M(Y)$ is degenerate
either. Let ${\Cal K}\cap T$ be a non-degenerate polyhedron. Then
the polyhedron $\M(Y)\cap C_{\Cal K}=\M(Y)\cap C_{{\Cal K}\cap T}$ is
bounded by a finite set of hyperplanes of codimension one faces of
$C_{{\Cal K}\cap T}$ which all contain the point $\Pa$ and by
hyperplanes $\Ha_E,\ E\in \Exc (Y)_1^\prime$, which have a non-empty
intersection with the cone $C_{{\Cal K}\cap T}$ over the compact
set ${\Cal K}\cap T$ of the horosphere $\E_\Pa$. Moreover,
$0<E\cdot P\le 1$ because
$-E\cdot K=E\cdot P+E\cdot \sum_i{\alpha_iF_i}=1$ and
$E\cdot P > 0$, $E\cdot F_i\ge 0$, $\alpha_i>0$.
This is a simple purely arithmetic
statement valid for any hyperbolic lattice (here, $\NS(Y)$),
any element $0\not= P\in \NS(Y)\otimes {\Bbb Q}$ with
$P^2=0$ and any compact set ${\Cal K}\cap T\subset \E_\Pa$,
$\Pa =\r^+P$, that the set
$$
\{e \in \NS(Y) \mid e^2=-1,\ 0< e\cdot P\le 1,\
\Ha_e \cap C_{{\Cal K}\cap T}\not=\emptyset\}
$$
is finite. (One should remark that hyperplanes $\Ha_e$ with
fixed $\lambda= e\cdot P>0$ are touching some horosphere
$\E_{\Pa, R(\lambda)}$, and the set of tangent points is compact in
$\La(Y)$
if additionally $\Ha_e \cap C_{{\Cal K}\cap T}\not=\emptyset$.
See the proof of a similar Statement below.)
Thus, the polyhedron $\M(Y)\cap C_{{\Cal K}\cap T}$ is bounded by a finite
set of hyperplanes in $\La (Y)$. If this polyhedron is not elliptic,
there exists an open subset
$U\subset \La(Y)_\infty \cap \M(Y)\cap C_{{\Cal K}\cap T}$
which does not contain $\Pa$. Since
$\M(Y)-\{\Pa\} \subset P\nem _-(Y)$, we get a contradiction with Lemma 3.2.
Thus, the polyhedron $\M(Y)\cap C_{\Cal K}=\M(Y)\cap C_{{\Cal K}\cap T}$
is elliptic.

Case (d).  Let ${\Cal K}\subset \T$ be a compact elliptic polyhedron
in $\T$ and $C_{\Cal K}\cap \M(Y)$ is non-degenerate.

We know that
$\M(Y)$ is a locally finite polyhedron with acute angles and
$P(\M(Y))=\Exc (Y)$. The subset
$\{F_1,...,F_n\}=\Exc (Y)_3\subset \Exc (Y)$ has a negative definite
Gram matrix. It follows (see Sect. 1.3) that
$$
\T=\bigcap_{F \in Exc (Y)_3}{\Ha_F}
$$
is a subspace of $\La$ of the codimension $n$
and $\gamma =\T \cap\M(Y)$ is
a face of $\M(Y)$ of the codimension $n$. Moreover,
$\M(Y)-\gamma \subset P\nem _-(Y)$, since $h\cdot F_i>0$ for
at least one $i$, $1\le i \le n$, if
$\r^+h\in \M(Y)-\gamma$, and
$-K \equiv \sum \alpha_iF_i$ where
$\alpha_i>0$ and $h\cdot F_i\ge 0$. Let $E\in \Exc(Y)_2$.
Since $E\notin \Exc (Y)_3$, $-K=\sum_i{\alpha_iF_i}$ where $\alpha_i>0$,
and $E\cdot K=0$, we have $E\cdot \Exc (Y)_3=0$. Thus, hyperplanes
$\Ha_E,\ E\in \Exc (Y)_2$ are orthogonal to $\T$. Since $\M(Y)$ is
locally-finite, it follows that $C_{\Cal K}\cap \M(Y)$ is bounded by
a finite set of hyperplanes $\Ha_E,\ E\in \Exc(Y)_3\cup \Exc(Y)_2$.
Let $E\in \Exc (Y)_1$. Then
$$
1=E\cdot (-K)=\sum_{i=1}^n{\alpha_i(E\cdot F_i)}.
$$
Since $\alpha_i$ are rational, $\alpha_i>0$ and $E\cdot F_i\ge 0$,
we have $0\le E\cdot F_i\le N$ for some constant $N$
depending from the rational numbers $\alpha_i$ ($N$ is
not greater than the least common multiple of denominators of
$\alpha_1,...,\alpha_n$). Since the lattice $\NS(Y)$ is
hyperbolic, we can prove purely arithmetically the following

\proclaim{Statement} The
set
$$
\{ e \in \NS(Y) \mid e^2=-1,\ 0\le e\cdot F_i\le N\
\text{for all\ } 1\le i \le n, \
\text{and}\ \Ha_e\cap C_{\Cal K}\not=\emptyset \}
$$
is finite.
\endproclaim

\demo{Proof}
Let $S$ be a hyperbolic lattice (i.e. a hyperbolic integral symmetric
bilinear form) and elements
$f_i\in S,\ 1\le i \le n,$ have a negative definite Gram matrix
$(f_i\cdot f_j)$. Let
$\T=\bigcap_{i=1}^n{\Ha_{f_i}}\subset \La(S)$ be the corresponding
subspace in $\La (S)$ and
${\Cal K}\subset \T$ a ball with a center $Q\in \T$ and a radius $R$.
We should prove that
the set
$$
\{ e \in S \mid e^2=-1,\ 0\le e\cdot f_i\le N\
\text{for all\ } 1\le i \le n, \
\text{and}\ \Ha_e\cap C_{\Cal K}\not=\emptyset \}
$$
is finite.

The Gram matrix of elements $e, f_1,...,f_n$ determines the configuration
of hyperplanes $\Ha_e, \Ha_{f_1},...,\Ha_{f_n}$ up to motions of the
Lobachevsky space $\La(S)$. By our conditions, the set of possible
Gram matrices of $e,f_1,...,f_n$ is finite. So, we can assume that
the configuration of the hyperplanes
$\Ha_e, \Ha_{f_1},...,\Ha_{f_n}$ is fixed up to motions of $\La (S)$.
Thus, we can assume that the configuration of the hyperplane
$\Ha_e$ and the subspace $\T$ is fixed up to motions of $\La(S)$.
This configuration is defined by either the angle or the distance
between $\Ha_e$ and $\T$.

Let us assume that the hyperplane $\Ha_e$ intersects the  subspace $\T$
in a finite or an infinite point.
Considering a bigger ball, we can suppose that the center $Q$ of the
ball is $Q=\r^+h$ where $h\in S$ and $h^2>0$.  Let us consider a
$2$-dimensional Lobachevsky plane with the same  curvature $(-1)$ and
a line $AO$ in the plane
where $A$ is a point at infinity  and $O$ is a finite point.
Let $B \in [A,O]$ and
$\rho (B,O)=R$. Let $BC$ be the perpendicular line  to the line $AO$
where $C$ is a point at infinity.  We put
$\rho(R)=\rho (AC, O)$ the distance
between the line  $AC$ and the point $O$. Obviously,
the constant  $\rho (R)$ only depends from $R$.  Elementary geometrical
considerations show that  $$ \rho (\Ha_e, Q) \le \rho (R).  $$ It
follows that $\vert (e\cdot h)/(e^2h^2)^{1/2}\vert \le \cosh \rho (R)$.
 Since $e^2=-1$ and $h$ is a fixed element, we get  $\vert e\cdot
h\vert<M$ for some constant $M$. Since $S$ is a hyperbolic  (integral)
lattice, the element $h\in S$ with $h^2>0$ is fixed, and $e^2=-1$ is fixed,
the set of these elements
$e\in S$ is finite.

Now we assume that the distance $\rho(\Ha_e, \T)=r>0$ is fixed.
Let us consider the set
$$
G(r)=\{X\in \La(S) \mid \rho (X, \T)=r\}.
$$
This is the boundary of the strongly convex closed in $\La(S)$ set
$$
\widetilde{G(r)}=\{X\in \La(S) \mid \rho (X, \T)\le r\}.
$$
The set
$\overline{G(r)}=G(r)\cup \T_\infty$ is a closed hypersurface in
$\overline{\La(S)}$ tangent the
$\La(S)_\infty$ exactly in points of $\T_\infty$. Obviously,
$\Ha_e$ is a hyperplane touching $G(r)$ exactly in one point of $G(r)$.

We consider an orthogonal basis $\xi_0,\xi_1,...,\xi_m$ in
$S\otimes \r$ such that $\xi_0^2=1$ and $\xi_i^2=-1$ for $1\le i \le m$.
We can assume that $\T$ is orthogonal to
$\xi_{m-n+1},...,\xi_m$, and the point $Q = \r^+\xi_0$. Then
a point of $\overline{\La(S)}$ is uniquely defined by a ray
$\r^+v$ where $v=\xi_0+x_1\xi_1+ \cdots +x_m\xi_m$, and one can consider
$(x_1,...,x_m),\ x_i\in \r$, as coordinates in $\overline{\La(S)}$.

In these coordinates, the $\overline{\La(S)}$ is the ball
$$
\overline{\La(S)}:\ \ x_1^2+\cdots+x_m^2\le 1,
$$
the subspace $\T$ is defined by the system of equations in $\La(S)$
$$
\T:\ \ x_{m-n+1}=...=x_m=0;
$$
the hypersurface $G(r)\cup \T_\infty$ is the ellipsoid
$$
G(r):\ \  x_1^2+\cdots +x_{m-n}^2+
x_{m-n+1}^2/\lambda (r)+\cdots +x_m^2/\lambda (r)=1
$$
where $\lambda (r)<1$.
The cylinder $C_{\Cal K}$ is
defined by the inequality in $\overline{\La(S)}$
$$
C_{\Cal K} :\ \ x_1^2+\cdots +x_{m-n}^2\le \mu (R)
$$
where $\mu (R)<1$. For this model, a hyperplane in $\La(S)$
is a section of $\La (S)$ by an affine linear hyperplane
in the coordinates
$x_1,...,x_m$.
Thus, hyperplanes $\Ha_e$ are tangent
to the ellipsoid $G(r)$ hyperplanes which intersect
the cylinder $C_{\Cal K}$ in $\overline{\La (S)}$.
Using the above coordinate description, one can easily find that the set of
common
points of the hyperplanes $\Ha_e$ and $G(r)$ is contained in
the compact subset $A\subset G(r)$ which is defined by the condition:
$x\in A$ if and only if $x\in G(r)$ and
the tangent to $G(r)$ at the point $x$ hyperplane
has a common point with
the cylinder $C_{\Cal K}\cap \overline{\La(S)}$. Let $A_\epsilon$ be the open
$\epsilon$-neighborhood of $A$ in $\La(S)$.
Then all hyperplanes $\Ha_e$
intersect the open subset $A_\epsilon \subset \La(S)$ with the compact closure
$\overline{A_\epsilon}$ in $\La(S)$.
Since the set of hyperplanes $\Ha_e$
with $e \in S$ and $e^2=-1$ is locally finite in $\La(S)$,
it follows, that the set of hyperplanes $\Ha_e$ under consideration
is finite.
$\blacksquare$
\enddemo

We continue the proof of Proposition 3.1.
By Statement above, the polyhedron $C_{\Cal K}\cap \M(Y)$ is
bounded by the finite set of hyperplanes $\Ha_e,\ e\in \Exc(Y)_1$.
By considerations before the statement,
it follows that the polyhedron $C_{\Cal K}\cap \M(Y)$ is bounded by
the finite set of hyperplanes.
If $C_{\Cal K}\cap \M(Y)$ is not elliptic, there is a non-empty open
subset $D\subset C_{\Cal K}\cap \M(Y)\cap \La_\infty$. By construction,
the polyhedron $C_{\Cal K}\cap \M(Y)$ is elliptic in a neighborhood of
$\overline \T$. Thus, we can assume that $D\cap \overline \T=\emptyset$.
Since $\M(Y)-\M(Y)\cap \overline \T \subset P\nem_- (Y)$, it follows
that
$D\subset \M(Y) \cap P\nem _-(Y) \cap \La (Y)_\infty$.
This contradicts to Lemma 3.2. Thus, $\M(Y)\cap C_{\Cal K}$ is elliptic.
$\blacksquare$
\enddemo

We can apply Proposition 3.1, to describe Mori cone $\nem (Y)$.
We recall that a ray $\r^+\delta \subset \nem (Y),
0\not=\delta \in \nem (Y)$, is called extremal if
$\delta=\delta_1+\delta_2$ where $0\not=\delta_1 \in \nem (Y)$ and
$0\not=\delta_2 \in \nem (Y)$ implies that
$\delta_1, \delta_2 \in \r^+\delta$. Obviously, $\nem (Y)$ is
generated by extremal rays: for any $0\not= x \in \nem (Y)$ there are
extremal rays $\r^+\delta_1,...,\r^+\delta_k$ of $\nem (Y)$ such that
$x=\delta_1+\cdots +\delta_k$.

An extremal ray $\r^+\delta$ of $\nem (Y)$ gives rise for
the dual polyhedron $\r^+\M(Y)$ the
so called supporting half-space $\Ha_\delta^+$.
That is the half-space
$$
\Ha_\delta^+=\{\r^+ x \subset \NS(Y)\otimes \r\mid
x\not= 0\ \text{and\ } \delta \cdot x\ge 0\}
$$
which contains $\M(Y)$ and there don't exist non-zero
$\delta_1, \delta_2 \in NS(Y)\otimes \r$ such that
$\Ha_{\delta_1}^+\not=\Ha_\delta^+,\
\Ha_{\delta_2}^+ \not=\Ha_\delta^+$,
$\M(Y)\subset \Ha_{\delta_1}^+\cap \Ha_{\delta_2}^+$ and
$\Ha_{\delta_1}^+\cap \Ha_{\delta_2}^+ \subset \Ha_\delta^+$.
Evidently, $\r^+\delta$ is an extremal ray of $\nem (Y)$ if and only if
$\Ha^+_\delta$ is a supporting half-space of $\M(Y)$. Thus,
the description of extremal rays of $\nem (Y)$ is equivalent
to the description of supporting half-spaces of $\M(Y)$.

In particular, since $\M(Y)\subset P\overline{V^+(Y)}$, we get
that $\delta^2\le 0$ for an extremal ray $\r^+\delta$ of
$\nem (Y)$. Using the description of extremal rays as supporting
half-spaces of $\M(Y)$ and Proposition 3.1, we get

\proclaim{Corollary 3.3} Let $Y$ be a non-singular projective
algebraic surface. Then Mori polyhedron $\nem (Y)$
is generated by extremal rays
$\r^+\delta$ where $\delta \in \NS(Y)\otimes \r$ and $\delta^2\le 0$.

Assume that the canonical class $K=K_Y\not\equiv 0$ and
the numerical anti-Kodaira dimension $\nu (-K,Y)\ge 0$. Let
$-K\equiv P+N,\ N=\sum_{i=1}^n{\alpha_iF_i},$ be Zariski
decomposition of $-K$. Then

(a) Extremal rays $\r^+\delta$ of $\nem (Y)$ where $\delta^2<0$ are
exactly extremal rays $\r^+E$ where $E$ is an exceptional curve.

(b) If $\nu (-K,Y)=2$ (i.e. $P^2>0$), there is no an extremal ray
$\r^+\delta$ with $\delta^2=0$; the set $\Exc (Y)$ is
finite and generates $\nem (Y)$.

(c) If $\nu (-K, Y)=1$ (i.e. $P^2=0$ but $P\not\equiv 0$),
there may exist the only extremal ray $\r^+c$ with $c^2=0$: this is
$\r^+P$. The ray $\r^+P$ is extremal if and only if $P\not=
\sum_{j=1}^{m}{a_jE_j}$ with $a_j>0$ and $E_j\in \Exc (Y)$. In
particular, $\nem (Y)$ is generated by $\Exc (Y)\cup \{P\}$.

(d) If $\nu (-K,Y)=0$ (i.e. $P\equiv 0$) and $\r^+c$ an extremal ray
of $\nem (Y)$ with $c^2=0$, then $\r^+c \in \M(Y)$ and
$c\cdot F_i=0$ for all $i=1,...,n$ (thus, $\r^+c \in \M(Y)\cap \T$
where $\T =\bigcap_{i=1}^n{\Ha_{F_i}}$).
Backwards, a ray $\r^+c \in \M(Y)\cap \T$ is extremal if and only if
$c \not=\sum_{j=1}^m{a_jE_j}$ where $a_j>0$ and $E_j\in \Exc (Y)$.
In particular, $\nem (Y)$ is generated by
$\Exc(Y)\cup \r^+(\M(Y)\cap \T \cap \La (Y)_\infty)$.
\endproclaim

\demo{Proof} We don't need this statement for the diagram method and
leave details of the proof to reader.
$\blacksquare$
\enddemo

It is the most important for us that by Proposition 3.1, we
can apply the diagram method (results of Sect. 1.4) to $Y$
if $-K=-K_Y$ is pseudo-effective. We only remark that if
$\nu (-K, Y)=1$, then $\M(Y)$ is parabolic relative to $\r^+P$ and
the element $e$ of Theorem 1.4.2 (or Theorem 1.4.2') has property
$e\cdot P>0$. It follows that $e$ is the class of an exceptional curve
$E$ of the first kind. If $\nu (-K,Y)=0$, when $\M(Y)$ is hyperbolic
relative to the subspace $\T =\bigcap_{i=1}^n{\Ha_{F_i}}$ of $\dim \T =
\dim \La(Y)-n$,
the subset $Q$ is a set
of classes of non-singular rational curves $E_1,...,E_m$ with square
$(-1)$ or $(-2)$, different from $F_1,...,F_n$,
and $m \le \dim \La -n$.
Theorems 1.4.3 and 1.4.3'give estimates for the number
$n$ of exceptional curves $F_1,...,F_n$.

Thus, applying Proposition 3.1 and results of Section 1.4, we get the
following Diagram method Theorems for surfaces
$Y$ with $\nu (-K, Y)\ge 0$.

\proclaim{Theorem 3.4} Let $Y$ be
a non-singular projective algebraic
surface over an algebraically closed field, and
numerical anti-Kodaira dimension
$\nu (-K, Y)=2$.
Assume that there are some constants $d, C_1, C_2$
such that the conditions (a) and (b) below hold:

(a)
For any Lanner subset
$L \subset \Exc (Y)$
$$
\text{diam}~\Gamma (L)\le d.
$$

(b) For any elliptic subset $\E\subset Exc (Y)$ such that
$\E$ has $\dim \NS (Y)-2$
elements, we have for the distance in
the graph $\Gamma (\E)$:
$$
\sharp \{ \{ E_1, E_2 \}\subset \E
\mid 1 \le \rho (E_1,E_2)\le d\}
\le C_1\sharp \E;
$$
and
$$
\sharp \{ \{E_1, E_2\} \subset \E
\mid d+1\le \rho (E_1,E_2) \le 2d+1\}
\le C_2 \sharp \E.
$$
Then $\dim \NS (Y) <96(C_1+C_2/3)+69$.
\endproclaim

\proclaim{Theorem 3.5} Let $Y$ be
a non-singular projective algebraic
surface over an algebraically closed field, and
numerical anti-Kodaira dimension
$\nu (-K, Y)=1$. Let $-K \equiv P+N$ be Zariski
decomposition where $P^2=0$ and $P\not\equiv 0$.
Assume that there are some constants $d, C_1, C_2$
such that the conditions (a) and (b) below hold:

(a)
For any Lanner subset
$L \subset \Exc (Y)$
$$
\text{diam}~\Gamma (L)\le d.
$$

(b) For any exceptional curve $E$ of the first kind such that
$E\cdot P>0$ and any
elliptic subset $\E \subset \Exc (Y)$
such that $\E$ has $\dim \NS (Y) -2$
elements and $E \in \E$, we have for the distance in
the graph $\Gamma (\E)$:
$$
\sharp \{ \{ E_1, E_2 \}\subset \E-\{E\}
\mid 1 \le \rho_E (E_1,E_2)\le d\}
\le C_1(\sharp \E-1);
$$
and
$$
\sharp \{ \{E_1, E_2\} \subset \E-\{E\}
\mid d+1\le \rho_E (E_1,E_2) \le 2d+1\}
\le C_2 (\sharp \E - 1).
$$
Then $\dim \NS(Y) <96(C_1+C_2/3)+70$.
\endproclaim

\proclaim{Theorem 3.6} Let $Y$ be a non-singular projective algebraic
surface over an algebraically closed field and numerical
anti-Kodaira dimension $\nu (-K, Y)=0$. Let
$-K\equiv \sum_{i=1}^n{\alpha_iF_i}$
be Zariski decomposition of $-K$ (i.e. all $\alpha_i>0$ and
the Gram matrix of $F_1,...,F_n$ is negative definite).
Assume that there are some constants $d, C_1, C_2$
such that the conditions (a) and (b) below hold:

(a)
For any Lanner subset
$L \subset \Exc (Y)$
$$
\text{diam}~\Gamma (L)\le d.
$$

(b) For any elliptic subset $Q\subset \Exc (Y)$ such that $Q$
contains only non-singular rational curves with square $(-1)$ or $(-2)$
which are different from $F_1,...,F_n$ above and
$\sharp Q \le \dim \NS(Y)-n-1$, and for any elliptic subset
$\E\subset \NS(Y)$ which contains
$\dim \NS(Y) -2$
elements and $Q \subset \E$, we have for the distance in
the graph $\Gamma (\E)$:
$$
\sharp \{ \{ E_1, E_2 \}\subset \E-Q
\mid 1 \le \rho_Q (E_1,E_2) \le d\}
\le C_1\sharp (\E-Q);
$$
and
$$
\sharp \{ \{E_1, E_2\} \subset \E-Q
\mid d+1 \le \rho_Q (E_1,E_2) \le 2d+1\}
\le C_2 \sharp (\E - Q).
$$
Then $n <96(C_1+C_2/3)+68$.
\endproclaim

\proclaim{Theorem 3.5'} Let $Y$ be
a non-singular projective algebraic
surface over an algebraically closed field and
numerical anti-Kodaira dimension
$\nu (-K, Y)=1$. Let $-K \equiv P+N$ be Zariski
decomposition where $P^2=0$ and $P\not\equiv 0$.
Assume that there are some constants $d, C_1, C_2$
such that the conditions (a) and (b) below hold:

(a)
For any Lanner subset
$L \subset \Exc (Y)$
$$
\text{diam}~\Gamma (L)\le d.
$$

(b) For any exceptional curve $E$ of the first kind such that
$E\cdot P>0$ and any
elliptic subset $\E \subset \Exc (Y)$
such that $\E$ has $\dim \NS (Y) -2$
elements and $E \in \E$, we have for the distance in
the graph $\Gamma (\E)$:
$$
\sharp \{ \{ E_1, E_2 \}\subset \E-\{E\}
\mid 1 \le \rho (E_1,E_2)\le d\}
\le C_1(\sharp \E-1);
$$
and
$$
\sharp \{ \{E_1, E_2\} \subset \E-\{E\}
\mid d+1\le \rho (E_1,E_2) \le 2d+1\}
\le C_2 (\sharp \E - 1).
$$
Then $\dim \NS(Y) <96(C_1+C_2/3)+70$.
\endproclaim

\proclaim{Theorem 3.6'} Let $Y$ be a non-singular projective algebraic
surface over an algebraically closed field and numerical
anti-Kodaira dimension $\nu (-K, Y)=0$. Let
$-K=\sum_{i=1}^n{\alpha_iF_i}$
be the Zariski decomposition of $-K$ (i.e. all $\alpha_i>0$ and
the Gram matrix of $F_1,...,F_n$ is negative definite).
Assume that there are some constants $d, C_1, C_2$
such that the conditions (a) and (b) below hold:

(a)
For any Lanner subset
$L \subset \Exc (Y)$
$$
\text{diam}~\Gamma (L)\le d.
$$

(b) For any elliptic subset $Q\subset \Exc (Y)$ such that $Q$
contains only non-singular rational curves with square $(-1)$ or $(-2)$
which are different from $F_1,...,F_n$ above and
$\sharp Q \le \dim \NS(Y)-n-1$, and for any elliptic subset
$\E\subset \NS(Y)$ which contains
$\dim \NS(Y) -2$
elements and $Q \subset \E$, we have for the distance in
the graph $\Gamma (\E)$:
$$
\sharp \{ \{ E_1, E_2 \}\subset \E-Q
\mid 1 \le \rho (E_1,E_2) \le d\}
\le C_1\sharp (\E-Q);
$$
and
$$
\sharp \{ \{E_1, E_2\} \subset \E-Q
\mid d+1 \le \rho (E_1,E_2) \le 2d+1\}
\le C_2 \sharp (\E - Q).
$$
Then $n <96(C_1+C_2/3)+68$.
\endproclaim

\remark{Remark 3.7} Probably, it is the most interesting applying
Theorems 3.4--3.6 and 3.5', 3.6' to normal projective algebraic
surfaces $Z$ with numerically effective anti-canonical class. We refer
to Mumford \cite{25} and Sakai \cite{28} on intersection theory of
Weil divisors on a normal surface.

Let $Z$ be a normal projective algebraic surface and
$\sigma : Y \to Z$ be a resolution of singularities of $Z$. Then
canonical classes are connected by the formula:
$$
K_Y \equiv \sigma ^\ast K_Z + \sum_{i=1}^n{\alpha_iF_i},
\tag4
$$
where $F_i$ are irreducible components of the exceptional divisor of
the resolution (i.e. $\sigma (F_i)$ is a singular point). One says that
$\sigma$ is minimal if any $F_i$ is not an exceptional curve of the first
kind. Then $\alpha_i\le 0$ for all $i$. This property is very important
for us. So, we name a resolution of singularities {\it almost minimal}
if in the formula (4) corresponding to this resolution all $\alpha_i \le 0$.

Let us assume that the surface $Z$ has $nef$ anti-canonical class $-K_Z$.
Then $(-\sigma^\ast K_Z)$ is $nef$, and for the almost minimal resolution
of singularities $\sigma$, the formula (4) gives the Zariski decomposition
of $-K_Y$ (if we take away zero summands in the $\sum$):
$$
-K_Y \equiv -\sigma^\ast K_Z +\sum_{i=1}^n{(-\alpha_i)F_i},
\tag5
$$
where $P=-\sigma^\ast K_Z$ is $nef$ (since
$-K_Z$ is $nef$) and
$N=\sum_{i=1}^n{(-\alpha_i)F_i}$. Here $(-\alpha_i)\ge 0$ since
$\sigma$ is almost minimal, the Gram matrix of $F_1,...,F_n$ is
negative definite by Mumford's Theorem (see \cite{25}) and
evidently $P$ is orthogonal to all $F_i$. Since $P^2=
(-\sigma^\ast K_Z)^2=(K_Z)^2$, we obtain
$$
\nu (-K_Y, Y)=
\cases
2, &\text{if $(K_Z)^2>0$};\\
1, &\text{if $(K_Z)^2=0$ and $K_Z\not\equiv 0$};\\
0, &\text{if $K_Z \equiv 0$};
\endcases
\tag6
$$
Applying Theorems 3.4--3.6 and 3.5', 3.6' to the surface $Y$,
we get a bound for $\dim \NS(Y)$ for elliptic and parabolic cases
(when $K_Z \not\equiv 0$). And we get a bound for
$n$ in the formula (5) if we take away all summands in (5) with $\alpha_i=0$.
In particular, for the minimal resolution of singularities of $Z$,
it is well-known that $\alpha_i =0$ if and only if $\sigma (F_i)$ is
a double rational point (Du Val singularity) of the type
${\Bbb A}_n, {\Bbb D}_n, {\Bbb E}_6, {\Bbb E}_7$ or ${\Bbb E}_8$.
Thus, for $K_Z\equiv 0$ we get a restriction on
singularities of $Z$ (in particular, on the number of non-Du Val
singularities).

\endremark

\enddocument

\end